\documentclass[12pt]{amsart}

\usepackage{graphicx,caption,subcaption}
\usepackage{tikz}        
\usepackage{amssymb}
\usepackage{epstopdf}
\usepackage[mathscr]{eucal}
\usepackage{amsmath,amssymb,amscd}
\usepackage{color}
\usepackage{epsfig}
\usepackage{hyperref}
\usepackage{cite}
\newtheorem{theorem}{Theorem}

\newtheorem{definition}{Definition}

\theoremstyle{definition}
\newtheorem{remark}{Remark}

\def \and{\mbox{and}}
\usepackage[top=4cm, bottom=4cm, left=3.5cm, right=3.5cm]{geometry}
\title{Linear stability of homogeneous and quasi-homogeneous N-body problem by symmetry groups}
\date{\today}

\begin{document}

\maketitle

\markboth{  }{  }
\author{\begin{center} 
Yingli Li ~~\footnote{School of Mathematics (Zhuhai), Sun Yat-sen University, Zhuhai, Guangdong, China \quad liyli26@mail2.sysu.edu.cn}\\
\end{center}}

\begin{abstract}
  Motivated by Xia-Zhou's recent work on applying symmetry groups to the N-body problem, 
  we will study relative equilibria of the equilateral triangle and the square configurations under $\alpha$-homogeneous and quasi-homogeneous potentials with this method.
 After linearizing the corresponding second order equations,
 with appropriate coordinate transformations, 
 we study the linear stability of the relative equilibria by decomposing each $2n\times 2n$ matrix into a series of $2\times 2$ matrices.  
 \end{abstract}
 

 \maketitle

 \section{{Introduction}}\label{sec:intro}
 In 2008, Xia\cite{Xia2008} provided insights on how the symmetry groups can be utilized to the stability of the N-body problem. 
 Recently, Xia-Zhou\cite{Xia2021} successfully applied this method to the linear stability of equal mass of three and four body problems.
 In this paper, we will study the linear stability of the corresponding relative equilibria $\alpha$-homogeneous and quasi-homogeneous potentials by systematically applying this method.

 The problem under $\alpha-$homogeneous potential is a natural generalisation of the Newtonian gravitational attraction, where $\alpha=1$.
Many properties of the Newtonian N-body problem (cf.\cite{Saari2005}\cite{Siegel1971})  have comparable properties in the homogeneous case. Furthermore, we will extend them to quasi-homogeneous N-body problem. 
 
 In Newtonian N-body problem, it is well-known that a relative equilibrium possesses very special configuration, which is called a central configuration. 
 It is hard to find central configurations as one needs to solve highly nonlinear algebraic equations\cite{Ferrario2008}.
 There are many open problems about central configurations. For example,
 relative equilibria of the three-body problem were solved by Euler and Lagrange,
 whereas for more than three body problem even the finiteness of central configurations is still an open problem.(cf.\cite{Albouy1996}\cite{Albouy2012}\cite{Hampton2006})
 
 In Barutello\cite{Vivina2014}, they use spectral flow \cite{Atiyah1975} to prove the following theorem.
   \begin{theorem}\cite{Hu2009}
     Let $\bar{x} $ be a critical point of the augmented potential and assume that it has even nullity. If the Morse index of $\bar{x}$ is odd, then the relative equilibrium corresponding to $\bar{x} $ is spectrally unstable.
   \end{theorem} 
 For $\alpha$-homogeneous potential, they demonstrated that the first eight eigenvalues  are always degenerate, if the relative equilibrium is non-degenerate, i.e.the remaining $4n-8$ eigenvalues are different from 0.
 In other words, the $4n-8$ eigenvalues are linked to the essential part of the dynamic.
 Take Lagrangian triangle with equal mass under  $\alpha$-homogeneous potential as an example,
 they concluded that the equilateral triangle is spectrally unstable for any $\alpha\in(0,2)$, this accords with the fact every regular polygon is linearly unstable in the gravitational case.(cf.\cite{Vivina2014}\cite{Moeckel1995})

 In Xia-zhou\cite{Xia2021}, they analysed the Newtonian N-body problem with certain symmetries, and the local stability of certain relative equilibria corresponding to symmetric central configurations.
 To be more specific, they considered two examples, the equilateral triangle with a fourth mass at the center and the square configurations with equal mass. 
 After linearizing the second order equation, they used group representation theory to decompose the hessian matrix $H$ into several lower order matrices in appropriate coordinates, leading to a series of lower order equations.
 As a result, eigenvalue calculations are effectively simplified and accurate.
 
 In Santoprete\cite{Santoprete2006}, they concisely proved that linear stability is achieved only when $a<2$, which is based on the idea of Moeckel(1992)\cite{Moeckel1992}.
 And they concluded that relative equilibria of the equilateral triangle configuration can be spectrally stable for $0<a<2$, when one of the masses is significantly larger than the others.
 Thus they inferred that if $a>0$, the case, where $m_{1}=m_{2}=m_{3}$, is always unstable.
 
 By applying the group representation theory, we systematically examine the linear stability of the equilateral triangle and the square configurations under $\alpha$-homogeneous potentials of relative equilibria with equal masses,
 and we obtain that for $\alpha>0$, these relative equilibria are unstable.

For configurations with potential functions $U=U_{1}+U_{2}$, where $U_{1}$ and $U_{2}$ are homogeneous functions of degree $-a$ and $-b$, respectively, with $0< a <b$. This type of function U is known as quasi-homogeneous potential.
Since the Manev potential, defined as $U(r)\sim \frac{1}{r} +\frac{1}{r^2} $, where $r$ is the relative distance between two mass points, is an alternative to the Newtonian potential that has a wide variety of applications in astrophysics and mechanics\cite{Manev2009}.
And the Schwarzschild potential, defined as $U(r)\sim \frac{1}{r} +\frac{1}{r^3} $, models astrophysical and stellar dynamic systems in a classical context \cite{Arredondo2014}.
 Therefore we choose these potentials to analyse the linear stability of the relative equilibria of the equilateral triangle and the square configurations with these potentials,
 in order to see if a similar conclusion can be reached. 
 As a result, we show that the linear stability of the equilateral triangle and the square configurations under the Manev and the Schwarzschild potentials of relative equilibria with equal mass are unstable.

 The paper is organized as follows. 
 In section \ref{sec:intro}, we give the background of the N-body problem and the linear stability of relative equilibria of three and four body problems. In addition, we compare the methods and results of various researches. 
 In section \ref{sec:prelim}, we review the basic concepts of central configuration, Hamiltonian of the N-body system and the group theory, all of which will be used in subsequent analysis.
 
 In section \ref{sec: the stab of homo}, using the group representation theory, we analyse the linear stability of relative equilibria in the N-body problem. We obtain a second-order equation from the Hamilton system. Furthermore, with appropriate coordinates, we linearize the second-order equation, greatly simplifying our calculations.
 In section \ref{sec: the stab of tria homo} and section \ref{sec: the stab of square homo}, we calculate the linear stability of relative equilibria of the equilateral triangle and the square configurations under homogeneous potential. With initial values and corresponding potentials, we get the Hessian matrices $H_{1}$ and $H_{2}$. Applying the representation theory, we get different trace equations, which help us find new coordinates. After linearizing the second-order equations, we can solve all the eigenvalues. By analysing them, we conclude that these configurations are unstable.
 
 From section \ref{sec: the stab of tria Manev} to section \ref{sec: the stab of square Schw}, we study the linear stability of the equilateral triangle and the square configurations of the Manev and Schwarzschild potentials, respectively, in detail. Unlike former homogeneous potential cases, we can not use Euler formula directly, but the original definition of central configurations to get the frequency $\omega^2$. Then we apply similar methods to calculate the Hessian matrices $H_{i},i=3,4,5,6$. By analysing all the eigenvalues, we find these configurations are unstable.
 \section{{Preliminaries}}\label{sec:prelim}
 In this section, we briefly review central configuration and the group theory that we will use in subsequent sections.
 
 \subsection{Central configuration}\label{sec:cen}
 We consider the equation
 \begin{equation}\label{equ:cen1}
   m_{i} \ddot{q}_{i}=\partial _{i}U =\alpha \sum_{1\le i<j\le n}^{} \frac{m_{i} m_{j}\left(q_{j}-q_{i}\right)}{\left\|q_{i}-q_{j}\right\|^{\alpha+2}} ,
   \end{equation}
   where $U = \sum_{1\le i<j\le n}^{} \frac{m_{i} m_{j}}{\left\|q_{i}-q_{j}\right\|^{\alpha}}$.

 Let the location of particles be $q_{i} \in \mathbb{R}^{2} $, $i=1,2,\dots,n$, let $q_{i}(t)=\phi(t) a_{i}$, where $a_{i}$ is a constant vector and $\phi(t)$ is a scalar vector. 
 A solution of this form is known as a homothetic solution.
 
 Then the equation is
 \begin{equation}\label{equ:cen2}
   \lambda_{\alpha} m_{i} a_{i}=\nabla U\left(a_{i}\right), \alpha>0,
 \end{equation}
 where $ \ddot{\phi}=-\frac{\lambda_{\alpha}}{\phi ^{\alpha+1}}$, $ \lambda_{\alpha}=\frac{\alpha U_{\alpha}\left(a_{i}\right)}{2 I\left(a_{i}\right)}$. 
 
\begin{definition}
  A central configuration is a configuration satisfies equation (\ref*{equ:cen2}).
\end{definition}

It is well-known that central configuration is the critical point of function $I^{\frac{\alpha}{2}}U $, where I satisfies
\begin{equation}
  I=\frac{1}{2} \sum_{i=1}^{n} m_{i}\left|q_{i}\right|^{2}.
  \end{equation}
The function $I$ is the moment of inertia of the N-body problem, the function $U$ is the $\alpha$-homogeneous self-potential function.

In particular, on the set $S=\left\{\vec{q}=\left(q_{1}, \ldots, q_{n}\right) \in R^{2 n} \mid I=1\right\}$, the central configuration is equivalent to the critical point of $U$.
    

\subsection{Hamiltonian of the N-body system}\label{sec:Hamilton2}

Let $p_{i}=m_{i} \dot{q}_{i}$ be the momentum, the hamiltonian of the N-body problem is
\begin{equation}
H=K-U,
\end{equation}
where 
\begin{equation}
K=\sum_{i=1}^{n} \frac{\left|p_{i}\right|^{2}}{2 m_{i}}.
\end{equation}

Thus the equation (\ref*{equ:cen1}) is equal to 
\begin{equation}
  \begin{cases}
    \dot{q}_{i}=\frac{\partial H}{\partial p_{i}}=\frac{p_{i}}{m_{i}},
     \\
    \dot{p}_{i}=-\frac{\partial H}{\partial q_{i}}=\frac{\partial U}{\partial q_{i}}.
    \end{cases}
  \end{equation}


\subsection{The group theory}\label{sec:group2}
\begin {definition}\label{def: representation of G}
Let $G$ be a group, a homomorphism $\mathscr{D}$ of $G$ into the group of $n\times n$ invertible matrices is a representation of $G$ of degree n.
\end {definition}

For any $A\in G, \mathscr{D}(A)$ is $n\times n$ non-singular matrix,
$\forall A_{1},A_{2}\in G$, we have
$\mathscr{D}(A_{1})\mathscr{D}(A_{2})=\mathscr{D}(A_{1}A_{2})$.

\begin{definition}\label{def:character of the representation}\cite{Steinberg2012}
  Let $\mathscr{D}$ be a representation of a group $G$ of degree $n$.
   $\forall$ $A\in G$, let $\chi(A)=Tr(\mathscr{D}(A))$, 
   thus $\chi$ defines a complex valued function on group $G$, 
   this function is called the character of the representation.
   
   In particular, representations with the same character are necessarily equivalent.
\end{definition}

\begin{definition}\label{def:inner product}\cite{Xia2008}
  Character of irreducible representation form an orthonormal basis under suitable defined inner product. Let $\chi_{i}$ and $\chi_{j}$ be characters of representations of a finite group  $G$.  We define the inner product by the formula
  \begin{equation}
    \left(\chi_{i}, \chi_{j}\right)=\frac{1}{|G|} \sum_{A \in G} 
    \overline{\chi_{i}(A)} \chi_{j}(A),
  \end{equation}
\end{definition}

where $G$ is a finite group and $|G|$ is the number of elements in $G$.

\begin{remark}\cite{Xia2008}
  This definition can be extended to infinite groups.
\end{remark}

With the above inner product of group representation, we have that if $\chi_{1}, \chi_{2}, \cdots, \chi_{k}$
    are characters of distinct irreducible representations, then
    \begin{equation}
      \left(\chi_{i}, \chi_{j}\right)=\delta_{i j}, i, j=1,2, \cdots, k.
    \end{equation}

    \begin{theorem}\label{thm:character calculation }\cite{Xia2008}
      Let $\mathscr{D}_{1}, \mathscr{D}_{2}, \cdots, \mathscr{D}_{h}$ be a set of irreducible representations of finite group $G$ and $h$ is the number of conjugacy classes of $G$.
      Let $\mathscr{D}_{1}, \mathscr{D}_{2}, \cdots, \mathscr{D}_{h}$ be the characters of these representations, with character $\chi$, any representation $\mathscr{D}$ of $G$ satisfies  
      \begin{equation}
        \mathscr{D}\sim n_{1} \mathscr{D}_{1} \oplus n_{2} 
        \mathscr{D}_{2} \oplus \cdots \oplus n_{h} \mathscr{D}_{h}.
      \end{equation}
    \end{theorem}

  \begin{definition}\label{def:invariant group}\cite{Xia2021}
    Let G be the finite group, A is the element of G, A acts on some vector space as a linear operator.
    We assume H is another linear operator defined on the same vector space, then H is invariant under the group G, if 
    \begin{equation*}
      \forall A\in G, HA=AH.
    \end{equation*}
  \end{definition}

    We consider the eigenspace $V_{\lambda}$ of H, $\forall \phi \in V_{\lambda}$, we obtain
    \begin{equation*}
      H \phi=\lambda \phi \Rightarrow A H \phi=\lambda A \phi \Rightarrow H(A \phi)=\lambda A \phi.
      \end{equation*}   
  
  $A\phi$ is an eigenvector for $H$ with the same eigenvalue.
  Denote $\phi_{1},\phi_{2},\cdots,\phi_{k}$ as a basis of the eigenspace $V_{\lambda}$, then
  \begin{equation*}
    A \phi_{j}=\sum_{i=1}^{k} \mathscr{D}_{i j}(A) \phi_{i}, \quad j=1,2, \ldots, k.
   \end{equation*}
  Let $\mathscr{D}(A)$ be a matrix that constitutes by  $\mathscr{D}_{i j}(A)$, for $A,B\in G$, we have
  \begin{equation}
    \mathscr{D}(B A) =\mathscr{D}(B) \mathscr{D}(A).
    \end{equation}
  So $\mathscr{D}$ is a group representation of $G$.\cite{Xia2021}
  
  We assume that the $n\times n$ matrix $H$ is the Hessian matrix of a smooth function of $n$ variables, and we choose n independent eigenvectors of $H$ to form the basis of  $ \mathbb{R}^n$.
  There is an invertible matrix P, fulfilling $\forall A\in G,P\mathscr{D}(A)P^{-1}=\mathscr{D}^\prime(A)$ and $PHP^{-1}=H^{\prime}$, where $H^{\prime}$ is formed by 
  $$
  H^{\prime}=\left(\begin{array}{cccc}
    \lambda_{1} & 0 & \ldots & 0 \\
    0 & \lambda_{2} & \ldots & 0 \\
    0 & 0 & \ldots & 0 \\
    \vdots & \vdots & \ddots & \vdots \\
    0 & 0 & \ldots & \lambda_{n}
    \end{array}\right).
  $$
  
  In the new coordinate system, the group representation is
  $$
  \mathscr{D}^{\prime}(A)=\left(\begin{array}{cccc}
    \mathscr{D}_{1}(A) & 0 & \ldots & 0 \\
    0 & \mathscr{D}_{2}(A) & \ldots & 0 \\
    0 & 0 & \ldots & 0 \\
    \vdots & \vdots & \ddots & \vdots \\
    0 & 0 & \ldots & \mathscr{D}_{k}(A)
    \end{array}\right).
  $$
  
  In addition, it satisfies
  \begin{equation}
    \operatorname{deg}\left(\mathscr{D}_{1}(A)\right)+\operatorname{deg}\left(\mathscr{D}_{2}(A)\right)+\cdots+\operatorname{deg}\left(\mathscr{D}_{k}(A)\right)=n. 
  \end{equation}
  
  Thus we obtain
  $$
  H\mathscr{D}(A)=H^{\prime}\mathscr{D}^{\prime}(A)=\left(\begin{array}{cccc}
    \lambda_{1} \mathscr{D}_{1}(A) & 0 & \ldots & 0 \\
    0 & \lambda_{2} \mathscr{D}_{2}(A) & \ldots & 0 \\
    0 & 0 & \ldots & 0 \\
    \vdots & \vdots & \ddots & \vdots \\
    0 & 0 & \ldots & \lambda_{k} \mathscr{D}_{k}(A)
    \end{array}\right).
  $$
  
  The group representation $\mathscr{D}$ is equivalent to $ \mathscr{D}_{1} \oplus \mathscr{D}_{2} \oplus \cdots \oplus \mathscr{D}_{k}$.\cite{Xia2021}

  \section{the stability of relative equilibria of the N-body problem}
  \label{sec:the stab of homogeneous potentials}
  
  \subsection{The stability of relative equilibrium of the N-body problem}\label{sec: the stab of homo}
  
  In this section, we analyse the stability of relative equilibria of the N-body problem in detail.
  
  When n particles move in concentric circles with the same angular velocity, then there is a series of special solutions known as relative equilibria.
  Let the center of concentric circle be at the origin, then
  \begin{equation}
  \begin{cases}
    q_{i}^{*}=\exp (-\omega J t) a_{i} ,
    \\
    p_{i}^{*}=-m_{i} \omega J \exp (-\omega J t) a_{i},
  \end{cases}
  \end{equation}
  where $a_{i}$ is a constant vector.
  
  While $J$ is the standard symplectic matrix
  $$
  J=\left(\begin{array}{cc}
    0 & 1 \\
    -1 & 0
    \end{array}\right).
  $$

 And the form of $ \exp (-\omega J t)$ is
  $$
    \exp (-\omega J t)=\left(\begin{array}{cc}
    \cos (\omega t) & -\sin (\omega t) \\
    \sin (\omega t) & \cos (\omega t)
    \end{array}\right),
  $$
  where $\omega$ is a positive constant satisfying the following equation. 
  \begin{equation}\label{equ:sec order equation}
  \omega^{2} a_{i}+\sum_{\substack{j=1 \\ j \neq i}}^{n} \frac{m_{j}\left(a_{j}-a_{i}\right)}{\left|a_{j}-a_{i}\right|^{\alpha}}=0.
  \end{equation}
  
  For the equation (\ref*{equ:sec order equation}), we perform linear transformation 
  \begin{equation}
    \begin{cases}
      q_{i}=\exp (-\omega J t)x_{i},
      \\
      p_{i}=\exp (-\omega J t)y_{i}.
    \end{cases}
  \end{equation}

  The equation is transformed as
  \begin{equation}
  \begin{cases}
    \dot{x_{i}}=\omega J x_{i}+\frac{y_{i}}{m_{i}},
    \\
    \dot{y_{i}}=\omega J y_{i}+\sum_{\substack{j=1 \\ j \neq i}}^{n} \frac{m_{i} m_{j}\left(x_{j}-x_{i}\right)}{\left|x_{j}-x_{i}\right|^{\alpha}}.
  \end{cases}
  \end{equation}
  
  As we know, the Hamilton system satisfies 
  \begin{equation}
    \begin{cases}
      \dot{x_{i}}=\frac{\partial H}{\partial y_{i}},
      \\
      \dot{y_{i}}=-\frac{\partial H}{\partial x_{i}}.
    \end{cases}
  \end{equation}
  
Thus with the same angular velocity, the Hamilton system is 
  \begin{equation}
  H=\sum_{i=1}^{n}\left(\frac{\left|y_{i}\right|^{2}}{2 m_{i}}-\omega x_{i}^{T} J y_{i}\right)-\sum_{1 \leq i<j \leq n} \frac{m_{i} m_{j}}{\left|x_{i}-x_{j}\right|^{\alpha}}.
  \end{equation}
  
 Both sides of this equation eliminate $y_{i}$ at the same time, which is transformed into the  second-order equation
  \begin{equation}
  \ddot{x}_{i}=2 \omega J \dot{x_{i}}+\omega^{2} x_{i}+\frac{1}{m_{i}} \sum_{\substack{j=1 \\ j \neq i}}^{n} \frac{m_{i} m_{j}\left(x_{j}-x_{i}\right)}{\left|x_{j}-x_{i}\right|^{\alpha}}.
  \end{equation}
  
  In order to simplify our calculation, we assume that the masses of all the particles are 1. By analysing this equation, we find that the first and second terms on the righ-hand side of the equation are linear, 
thus we only need to linearize the last term, which is related to potential $U$.
  This term is equivalent to the gradient of the potential, hence we mark it $\nabla_{x_{i}} U(x_{i})$.
  
Therefore the linearized second-order equation is
  \begin{equation}\label{equ19:linearized second equation}
    \ddot{x}_{i}=2 \omega J \dot{x_{i}}+\omega^{2} x_{i}+\nabla_{x_{i}} U, x\in \mathbb{R}^2.
  \end{equation}

    Let $X=\left(\begin{array}{c}
      x_{1} \\
      \vdots \\
      x_{n}
      \end{array}\right) \in \mathbb{R}^{2 n}$, substituting this into the equation (\ref*{equ19:linearized second equation}), we get  
      \begin{equation}
        \ddot{X}=2 \omega J \dot{X}+\omega^{2} X+D^{2} U\left(z_{0}\right) X,
      \end{equation}
  where $z_{0}$ is the corresponding central configuration.
   
  We find that if $\lambda$ is a double eigenvalue we can choose appropriate eigenvectors of $\lambda$ satisfying
      $J(V_{1},V_{2})=(V_{1},V_{2})J$.
     If $\lambda$ is a simple eigenvalue, we can combine with the other simple eigenvalue $\lambda_{i}$ so that
      $J(V_{1},V_{2})=(V_{1},V_{2})J.$ 
  
    If the eigenspace $E_{\lambda_{i}}=span\left \{ V_{1}, V_{2} \right \} $ of $\lambda_{i}$ satisfies
    $J(V_{1},V_{2})=(V_{1},V_{2})J$, then on each eigenspace, the linearized equation under this new basis is decomposed to
    \begin{equation}
      \ddot{x}=2 \omega J \dot{x}+\left(\omega^{2}I_{2}+\Lambda _{i}\right) x,
    \end{equation}
    where $\Lambda _{i}=\begin{pmatrix}
      \lambda _{i} & 0\\
      0 &\lambda _{j}
     \end{pmatrix}$. If $\lambda$ is a double eigenvalue, the index satisfies $i=j$.
  
  \subsection{The stability of relative equilibria of equilateral triangle problem with equal mass}\label{sec: the stab of tria homo}
  
  Let $ z \in \mathbb{R}^{6}$ satisfy $z=\left(x_{1}, y_{1}, x_{2}, y_{2}, x_{3}, y_{3}\right)^{T}$,
  we choose the vertices of the equilateral triangle $$z_{0}=(1,0,-\frac{1}{2},\frac{\sqrt{3}}{2},-\frac{1}{2},-\frac{\sqrt{3}}{2})^{T}$$
  as a critical point.
  
  Then we assume that the corresponding relative equilibria are
  
  $$
  \begin{array}{c}
    \mathbf{q}_{1}=\left(\begin{array}{c}
    \cos \omega t \\
    \sin \omega t
    \end{array}\right), \quad 
    \mathbf{q}_{2}=\left(\begin{array}{c}
    \cos \left(\omega t+\frac{2 \pi}{3}\right) \\
    \sin \left(\omega t+\frac{2 \pi}{3}\right)
    \end{array}\right) , \quad 
    \mathbf{q}_{3}=\left(\begin{array}{c}
    \cos \left(\omega t+\frac{4 \pi}{3}\right) \\
    \sin \left(\omega t+\frac{4 \pi}{3}\right)
    \end{array}\right). 
    \end{array}
  $$

  In $\alpha$-homogeneous three-body system, the degree of $U$ is $\alpha$, the degree of $I$ is 2, thus we apply Euler formula to calculate
  \begin{equation}
  \omega^2=\frac{\alpha U(z_{0})}{2I(z_{0})}=3^{-\frac{\alpha}{2}}\alpha.
  \end{equation}
  
  Then we obtain relative equilibria with initial configuration 
  $$z_{0}=(1,0,-\frac{1}{2}, \frac{\sqrt{3}}{2},-\frac{1}{2},-\frac{\sqrt{3}}{2})^{T}.$$
  
  We write the Hessian matrix $H_{1}:=\partial^{2} U(z) / \partial z^{2}$ at $z_{0}$ as 

  $H_{1}=
        \begin{pmatrix}
          \begin{smallmatrix}
    \frac{3^{-\frac{\alpha}{2}} \alpha(3+2 \alpha)}{6} & 0 & -\frac{3^{-\frac{\alpha}{2}} \alpha(3+2 \alpha)}{12}& \frac{3^{-\frac{\alpha}{2}} \alpha(2+\alpha)}{4\sqrt{3}} & -\frac{3^{-\frac{\alpha}{2}} \alpha(3+2 \alpha)}{12} & -\frac{3^{-\frac{\alpha}{2}} \alpha(2+\alpha)}{4\sqrt{3}} \\
    0 & \frac{3^{-\frac{\alpha}{2}} \alpha(-2+\alpha)}{6} & \frac{3^{-\frac{\alpha}{2}} \alpha(2+\alpha)}{4\sqrt{3}}& -\frac{3^{-\frac{\alpha}{2}} \alpha(-2+\alpha)}{12} & -\frac{3^{-\frac{\alpha}{2}} \alpha(2+\alpha)}{4\sqrt{3}} & -\frac{3^{-\frac{\alpha}{2}} \alpha(-2+\alpha)}{12} \\
    -\frac{3^{-\frac{\alpha}{2}} \alpha(3+2 \alpha)}{12}  & \frac{3^{-\frac{\alpha}{2}} \alpha(2+\alpha)}{4\sqrt{3}} & \frac{3^{-\frac{\alpha}{2}} \alpha(-2+3\alpha)}{12} & -\frac{3^{-\frac{\alpha}{2}} \alpha(2+\alpha)}{4\sqrt{3}} & \frac{3^{-\frac{\alpha}{2}} \alpha}{3} & 0 \\
    \frac{3^{-\frac{\alpha}{2}} \alpha(2+\alpha)}{4\sqrt{3}} & -\frac{3^{-\frac{\alpha}{2}} \alpha(-2+ \alpha)}{12}  & -\frac{3^{-\frac{\alpha}{2}} \alpha(2+\alpha)}{4\sqrt{3}} & \frac{3^{-\frac{\alpha}{2}} \alpha(2+5\alpha)}{12} & 0 & -\frac{3^{-\frac{\alpha}{2}} \alpha(1+\alpha)}{3} \\
    -\frac{3^{-\frac{\alpha}{2}} \alpha(3+2 \alpha)}{12}  & -\frac{3^{-\frac{\alpha}{2}} \alpha(2+\alpha)}{4\sqrt{3}} & \frac{3^{-\frac{\alpha}{2}} \alpha}{3}& 0 & \frac{3^{-\frac{\alpha}{2}} \alpha(-2+3\alpha)}{12}   & \frac{3^{-\frac{\alpha}{2}} \alpha(2+\alpha)}{4\sqrt{3}} \\
    -\frac{3^{-\frac{\alpha}{2}} \alpha(2+\alpha)}{4\sqrt{3}} & -\frac{3^{-\frac{\alpha}{2}} \alpha(-2+\alpha)}{12}  & 0 & -\frac{3^{-\frac{\alpha}{2}} \alpha(1+\alpha)}{3}  & \frac{3^{-\frac{\alpha}{2}} \alpha(2+\alpha)}{4\sqrt{3}} & \frac{3^{-\frac{\alpha}{2}} \alpha(5+2 \alpha)}{12} 
      \end{smallmatrix}
    \end{pmatrix}.
  $
   Clearly, $H_{1}$ is invariant under the group action $S_{3}$. Therefore we use the group representation theory to study the related eigenvalues and eigenvectors of $H_{1}$.
  
  Elements of the symmetric group $S_{3} $ are $3\times 3$ matrices as follows.
   
   $\, I=\begin{pmatrix}
       1 & 0 &0 \\
       0 & 1 &0 \\
       0 & 0 &1
      \end{pmatrix}$,\,\,\,\,\,\,\,
   $T=\begin{pmatrix}
       1 & 0 &0\\
       0 & 0 &1 \\
       0 & 1 &0
      \end{pmatrix}$,\,\,\,\,\,\,\,
     $R=\begin{pmatrix}
       0 & 0 &1\\
       1 & 0 &0 \\
       0 & 1 &0
      \end{pmatrix}$,
   
     $ R^{2} =\begin{pmatrix}
       0 & 1 &0\\
       0 & 0 &1\\
       1 & 0 &0
      \end{pmatrix}$,\
   $ TR =\begin{pmatrix}
       0 & 0 &1\\
       0 & 1 &0\\
       1 & 0 &0
      \end{pmatrix}$,
      $ TR^{2} =\begin{pmatrix}
       0 & 1 &0\\
       1 & 0 &0\\
       0 & 0 &1
      \end{pmatrix}$.
   
   Let $\mathscr{D}_{1}(A)=I$, $\mathscr{D}_{2}(A)=det(A)$, we deduce that $\mathscr{D}_{1}(A)$ and $\mathscr{D}_{2}(A)$
   are both of degree 1.
   
   And for the third group representation $\mathscr{D}_{3}(A)$ of $S_{3} $, its degree is 2. 
   
   Thus we define $\mathscr{D}_{3}(A)$ as follows.
   
   $\, \mathscr{D}_{3} (I)=\begin{pmatrix}
       1& 0\\
       0& 1
     \end{pmatrix}$,\,\,\,\
     $\mathscr{D}_{3} (T)=\begin{pmatrix}
       1& 0\\
       0& -1
     \end{pmatrix}$,
     \\
     $\mathscr{D}_{3} (R)=\begin{pmatrix}
       -\frac{1}{2}& -\frac{\sqrt{3}}{2}\\
       \frac{\sqrt{3}}{2}& -\frac{1}{2}
     \end{pmatrix}$,
    $ \mathscr{D}_{3} (R^{2})=\begin{pmatrix}
       -\frac{1}{2}& \frac{\sqrt{3}}{2}\\
       -\frac{\sqrt{3}}{2}& -\frac{1}{2}
     \end{pmatrix}$,
     \\
     $\mathscr{D}_{3} (TR)=\begin{pmatrix}
       -\frac{1}{2}& -\frac{\sqrt{3}}{2}\\
       -\frac{\sqrt{3}}{2}& \frac{1}{2}
     \end{pmatrix}$,
     $\mathscr{D}_{3} (TR^{2})=\begin{pmatrix}
       -\frac{1}{2}& \frac{\sqrt{3}}{2}\\
       \frac{\sqrt{3}}{2}& \frac{1}{2}
     \end{pmatrix}$.
   
   The group representation $\mathscr{D}_{3}(A)$ establishes an isomorphism between $S_{3} $ and $\mathscr{D}_{3}$.
   
   The character table of $S_{3} $  is
   
       $$
       \begin{array}{|l|l|l|l|l|}
           \hline A / \chi(A) & \chi_{1} & \chi_{2} & \chi_{3} & \chi_{4} \\
           \hline I & 1 & 1 & 2 & 3 \\
           \hline R & 1 & 1 & -1 & 0 \\
           \hline R^{2} & 1 & 1 & -1 & 0 \\
           \hline T & 1 & -1 & 0 & 1 \\
           \hline T R & 1 & -1 & 0 & 1 \\
           \hline T R^{2} & 1 & -1 & 0 & 1 \\
           \hline
           \end{array}
       $$
   
   For every element $A$ of group $S_{3} $, we apply linear transformation $\mathscr{D}(A)$ to express the effect of $A$ on $H_{1}$, 
   thus the transformation can be stated as $\mathscr{D}(A)H_{1}$.
   
   Since $H_{1}$ is invariant under the action of $S_{3}$, we obtain $$\forall A \in S_{3}, \mathscr{D}(A)^{-1} H_{1} \mathscr{D}(A)=H_{1}^{\prime}.$$
   
   Then for the group representation $\mathscr{D}$ of group $S_{3}$, the degree is 6.
   
   We list elements of the group representation $\mathscr{D}$ of group $S_{3}$, 
   where the matrix $R$ represents the entire configuration rotating $\frac{2\pi}{3}$ clockwise.

   $$
   \mathscr{D}(R)=\begin{pmatrix}
       0 & 0 & 0 & 0 & -\frac{1}{2} & -\frac{\sqrt{3}}{2} \\
       0 & 0 & 0 & 0 & \frac{\sqrt{3}}{2} & -\frac{1}{2} \\
       -\frac{1}{2} & -\frac{\sqrt{3}}{2} & 0 & 0 & 0 & 0 \\
       \frac{\sqrt{3}}{2} & -\frac{1}{2} & 0 & 0 & 0 & 0 \\
       0 & 0 & -\frac{1}{2} & -\frac{\sqrt{3}}{2} & 0 & 0 \\
       0 & 0 & \frac{\sqrt{3}}{2} & -\frac{1}{2} & 0 & 0
   \end{pmatrix}.
   $$
   
   And the matrix $T$ is a reflection about x-axis.
   
   $$
   \mathscr{D}(T)=\begin{pmatrix}
       1 & 0 & 0 & 0 & 0 & 0 \\
   0 & -1 & 0 & 0 & 0 & 0 \\
   0 & 0 & 0 & 0 & 1 & 0 \\
   0 & 0 & 0 & 0 & 0 & -1 \\
   0 & 0 & 1 & 0 & 0 & 0 \\
   0 & 0 & 0 & -1 & 0 & 0
   \end{pmatrix}.
   $$
   
   Representations of other elements of $S_{3}$ are generated by $\mathscr{D}(R)$ and $\mathscr{D}(T)$.
   
   Let $\chi(I)$ be the character of $\mathscr{D}$, then the character table is
   $$
   \begin{array}{|l|l|l|l|l|l|}
       \hline A / \chi(A) & \chi & \chi_{1} & \chi_{2} & \chi_{3} & \chi_{4} \\
       \hline I & 6 & 1 & 1 & 2 & 3 \\
       \hline R & 0 & 1 & 1 & -1 & 0 \\
       \hline R^{2} & 0 & 1 & 1 & -1 & 0 \\
       \hline T & 0 & 1 & -1 & 0 & 1 \\
       \hline T R & 0 & 1 & -1 & 0 & 1 \\
       \hline TR^{2} &0 & 1 & -1 & 0 & 1 \\
       \hline
       \end{array}
   $$
   From this character table, we obtain $\chi=\chi_{1}+\chi_{2}+2 \chi_{3}$, which means 
   $\mathscr{D}$ is equivalent to $\mathscr{D}_{1} \oplus \mathscr{D}_{2} \oplus \mathscr{D}_{3} \oplus \mathscr{D}_{3}$.
   
   Now suppose the matrix $H_{1}$ is in its Jordan canonical form, with suitable coordinates. 
   Using  $S_{3}$ symmetries in $H_{1}$, we group all multiple eigenvalues into the Jordan canonical form.
   Rearranging the order of eigenvalues, we assume that with new coordinates, the transformation $\mathscr{D}(A)$ is represented as 
   $$
   \mathscr{D}^{\prime}(A)=\left(\begin{array}{cccc}
     \mathscr{D}_{1}(A) & 0 & 0 & 0 \\
     0 & \mathscr{D}_{2}(A) & 0 & 0 \\
     0 & 0 & \mathscr{D}_{3}(A) & 0 \\
     0 & 0 & 0 & \mathscr{D}_{3}(A)
     \end{array}\right).
   $$
   We write the Jordan canonical form $H_{1}^{\prime}$ of matrix $H_{1}$ as follow.
   $$
    H_{1}^{\prime}=\left(\begin{array}{cccccc}
     \lambda_{1} & 0 & 0 & 0 & 0 & 0 \\
     0 & \lambda_{2} & 0 & 0 & 0 & 0 \\
     0 & 0 & \lambda_{3} & 0 & 0 & 0 \\
     0 & 0 & 0 & \lambda_{3} & 0 & 0 \\
     0 & 0 & 0 & 0 & \lambda_{4} & 0 \\
     0 & 0 & 0 & 0 & 0 & \lambda_{4}
     \end{array}\right).
   $$
   Then
   $$
   H_{1}\mathscr{D}(A)=H_{1}^{\prime}\mathscr{D}^{\prime}(A)=\left(\begin{array}{cccc}
     \lambda_{1} \mathscr{D}_{1}(A) & 0 & 0 & 0 \\
     0 & \lambda_{2} \mathscr{D}_{2}(A) & 0 & 0 \\
     0 & 0 & \lambda_{3} \mathscr{D}_{3}(A) & 0 \\
     0 & 0 & 0 & \lambda_{4} \mathscr{D}_{3}(A)
     \end{array}\right).
   $$

  The trace of $H_{1}\mathscr{D}(A)$ is
  \begin{equation*}
  \operatorname{Tr}(H_{1} \mathscr{D}(A))=H^{\prime}_{1}\mathscr{D}^{\prime}(A)=
  \lambda_{1} \chi_{1}(A)+\lambda_{2} \chi_{2}(A)+\lambda_{3} \chi_{3}(A)+\lambda_{4} \chi_{3}(A).
  \end{equation*}

  For $A=I$, the corresponding trace of $H_{1}\mathscr{D}(A)$ is
  
  \begin{equation*}
  \operatorname {Tr} (H_{1})=\lambda_ {1} + \lambda_ {2} + 
   2\left (\lambda_ {3} + \lambda_ {4} \right) =  2 \times 3^{-\frac{\alpha}{2} }\alpha^{2}.
  \end{equation*}
  
  For $A=R$, the corresponding trace of $H_{1}\mathscr{D}(R)$ is
  \begin{equation*}
  \operatorname {Tr} 
  (H_{1}\mathscr{D}(R))=\lambda_ {1} + \lambda_ {2} -\left (\lambda_ {3} + \lambda_ {4} \right) = \frac{1}{2} \times 3^{-\frac{\alpha}{2} } \alpha^{2}.
  \end{equation*}
  
  For $A=T$, the corresponding trace of $H_{1}\mathscr{D}(T)$ is
  \begin{equation*}
    \operatorname {Tr} 
    (H_{1}\mathscr{D}(T))= \lambda_ {1} - \lambda_ {2} + 
   0\left (\lambda_ {3} + \lambda_ {4} \right) = 3^{-\frac{\alpha}{2} }  \alpha(2+\alpha).
  \end{equation*}
  
  We solve the above equations by combining them.
  $$\lambda_ {1}=3^{-\frac{\alpha}{2}} \alpha(1+\alpha), \lambda_ {2}=-3^{-\frac{\alpha}{2}} \alpha, \lambda_ {3}+\lambda_ {4}= \frac{1}{2} \times 3^{-\frac{\alpha}{2}} \alpha^{2}.
  $$
  
Since at least two eigenvalues are equal to 0, thus we assume $\lambda_ {3}=0$.
  
  Substituting this assumption to the characteristic polynomial, we obtain
  
  $$\lambda_ {4}= \frac{1}{2} \times 3^{-\frac{\alpha}{2}} \alpha^{2}.$$
  
  In conclusion, the eigenvalues of matrix $H_{1}$ are as follows.
  
  \begin{equation}
    \begin{cases}
  \lambda_ {1} =3^{-\frac{\alpha}{2}} \alpha(1+\alpha),
  \\
  \lambda_ {2}=-3^{-\frac{\alpha}{2}} \alpha,
  \\
  \lambda_ {3}=0,
  \\
  \lambda_ {4}= \frac{1}{2} \times 3^{-\frac{\alpha}{2}} \alpha^{2}.
    \end{cases}
  \end{equation}
  where $\lambda_{3}$ and $\lambda_{4}$ are double eigenvalues, and others are simple eigenvalues.
  
  The forms of $\Lambda_{i}(i=1,2,3)$ are obtained by combining the eigenvalues.
  $$
  \begin{array}{c}
    \Lambda_{1}=\left(\begin{array}{cc}
    \lambda_{1} & 0 \\
    0 & \lambda_{2}
    \end{array}\right), \quad \Lambda_{2}=\left(\begin{array}{cc}
    \lambda_{3} & 0 \\
    0 & \lambda_{3}
    \end{array}\right), \quad 
    \Lambda_{3}=\left(\begin{array}{cc}
    \lambda_{4} & 0 \\
    0 & \lambda_{4}
    \end{array}\right).
  \end{array}
  $$
  
 Applying the group theory, we choose a series of eigenvectors as follows.

 For $\lambda_{1}$, we choose the eigenvector $v_{1}=\left\{1,0,-\frac{1}{2}, \frac{\sqrt{3}}{2},-\frac{1}{2},-\frac{\sqrt{3}}{2}\right\}^{T} $;

 For $\lambda_{2}$, we choose the eigenvector $v_{2}=\left\{0,1,-\frac{\sqrt{3}}{2},-\frac{1}{2}, \frac{\sqrt{3}}{2},-\frac{1}{2}\right\}^{T}$;
  
 The eigenvalue $\lambda_{3}$ is a double characteristic root, we choose two linearly independent eigenvectors, which are 
$v_{3}=\left\{0,1,0,1,0,1\right\}^{T}, v_{4}=\left\{1,0,1,0,1,0\right\}^{T}$;
  
The eigenvalue $\lambda_{4} $ is also a double characteristic root, 
 we choose two linearly independent eigenvectors, which are
  $$v_{5}=\left\{0,1, \frac{\sqrt{3}}{2},-\frac{1}{2},-\frac{\sqrt{3}}{2},-\frac{1}{2}\right\}^{T}, 
 v_{6}=\left\{1,0, -\frac{1}{2},-\frac{\sqrt{3}}{2},\frac{1}{2},\frac{\sqrt{3}}{2}\right\}^{T}.
  $$
  
  Therefore we choose these eigenvectors as a basis which satisfy
  \begin{equation}
    \begin{split}
  J\begin{pmatrix}
    v_{1} &v_{2}
   \end{pmatrix}&=\begin{pmatrix}
    0 &1 &0 &0 &0 &0\\
    -1&0 &0 &0 &0 &0\\
     0&0 &0 &1 &0 &0 \\
     0&0 &-1&0 &0 &0 \\
     0&0 &0 &0 &0 &1\\
     0&0 &0 &0 &-1&0
   \end{pmatrix}
   \left(\begin{array}{cc}
   1 & 0 \\
   0 & 1 \\
   -\frac{1}{2} & -\frac{\sqrt{3}}{2} \\
   \frac{\sqrt{3}}{2} & -\frac{1}{2} \\
   -\frac{1}{2} & \frac{\sqrt{3}}{2} \\
   -\frac{\sqrt{3}}{2} & -\frac{1}{2}
   \end{array}\right)=
   \left(\begin{array}{cc}
   0 & 1 \\
   -1 & 0 \\
   \frac{\sqrt{3}}{2} & -\frac{1}{2} \\
   \frac{1}{2} & \frac{\sqrt{3}}{3} \\
   -\frac{\sqrt{3}}{2} & -\frac{1}{2} \\
   \frac{1}{2} & -\frac{\sqrt{3}}{2}
   \end{array}\right) \\
  & =
   \left(\begin{array}{cc}
   1 & 0 \\
   0 & 1 \\
   -\frac{1}{2} & -\frac{\sqrt{3}}{2} \\
   \frac{\sqrt{3}}{2} & -\frac{1}{2} \\
   -\frac{1}{2} & \frac{\sqrt{3}}{2} \\
   -\frac{\sqrt{3}}{2} & -\frac{1}{2}
   \end{array}\right)
   \begin{pmatrix}
    0 &1 &0 &0 &0 &0\\
    -1&0 &0 &0 &0 &0\\
     0&0 &0 &1 &0 &0 \\
     0&0 &-1&0 &0 &0 \\
     0&0 &0 &0 &0 &1\\
     0&0 &0 &0 &-1&0
   \end{pmatrix}=\begin{pmatrix}
    v_{1} &v_{2}
   \end{pmatrix}J.
  \end{split}
  \end{equation}
  
  We also find that for $\lambda_{3}=0$, we have
  $
    J \begin{pmatrix}
      v_{3} &v_{4}
     \end{pmatrix}=\begin{pmatrix}
      v_{3} &v_{4}
     \end{pmatrix}J
  $;
  
 For $\lambda_{4}= \frac{1}{2} \times 3^{-\frac{\alpha}{2}} \alpha^{2}$, we have
  $
    J \begin{pmatrix}
      v_{5} &v_{6}
     \end{pmatrix}=\begin{pmatrix}
      v_{5} &v_{6}
     \end{pmatrix}J
  $.

  
  In the new coordinate $E_{\lambda}=\left \{ V_{i}, V_{j}\right \}$,$1\le i,j\le 6$, let $y_{i}=(x_{i},\dot{x_{i}})$, then the linearized second order equations are
  $\dot{y_{i}}=B_{i}y_{i},\,i=1,2,3$,
  where $B_{i}$ satisfy
  \begin{equation}
  B_{i}=\left(\begin{array}{cc}
    0 & I_{2} \\
    \omega^{2} I_{2}+\frac{1}{m_{i}} \Lambda_{i} & 2 \omega J
    \end{array}\right), i=1,2,3.
  \end{equation}
  
  We calculate all eigenvalues of equations $\dot{y_{i}}=B_{i}y_{i},\,i=1,2,3$.
  
  Firstly, we calculate the matrix $B_{1}$.
  \begin{equation}
    \begin{aligned}
  B_{1} &=\left(\begin{array}{cccc}
    0 & 0 & 1 & 0 \\
    0 & 0 & 0 & 1 \\
    \omega^{2}+\lambda_{1} & 0 & 0 & 2 \omega \\
    0 & \omega^{2}+\lambda_{2} & -2 \omega & 0
    \end{array}\right)\\
    &=\left(\begin{array}{cccc}
    0 & 0 & 1 & 0 \\
    0 & 0 & 0 & 1 \\
    3^{-\frac{\alpha}{2}}\alpha+3^{-\frac{\alpha}{2}} \alpha(1+\alpha) & 0 & 0 & 2 \sqrt{3^{-\frac{\alpha}{2}}\alpha} \\
    0 & 0 & -2 \sqrt{3^{-\frac{\alpha}{2}}\alpha} & 0
    \end{array}\right).
  \end{aligned}
  \end{equation}
  
By calculating eigenvalues of matrix $B_{1}$, we get 
  \begin{equation}
    \begin{cases}
  \lambda_ {1}^{\prime} =0,
  \\
  \lambda_ {2}^{\prime}=0,
  \\
  \lambda_ {3}^{\prime}=-\mathrm{i}\sqrt{ 3^{-\frac{\alpha}{2}}\alpha(2-\alpha)}=-\mathrm{i}\omega\sqrt{(2-\alpha)},
  \\
  \lambda_ {4}^{\prime}=\mathrm{i}\sqrt{ 3^{-\frac{\alpha}{2}}\alpha(2-\alpha)}=\mathrm{i}\omega\sqrt{(2-\alpha)};
    \end{cases}
  \end{equation}
  
  Similarly, we calculate $B_{2}$ to obtain
  \begin{equation}
  B_{2}=\left(\begin{array}{cccc}
    0 & 0 & 1 & 0 \\
    0 & 0 & 0 & 1 \\
    3^{-\frac{\alpha}{2}}\alpha & 0 & 0 & 2 \sqrt{3^{-\frac{\alpha}{2}}\alpha} \\
    0 & 3^{-\frac{\alpha}{2}}\alpha & -2 \sqrt{3^{-\frac{\alpha}{2}}\alpha} & 0
    \end{array}\right).
  \end{equation}
  
  By calculating eigenvalues of $B_{2}$, we get 
  \begin{equation}
    \begin{cases}
  \lambda_ {5}^{\prime} =-\mathrm{i}\sqrt{ 3^{-\frac{\alpha}{2}}\alpha}=-\mathrm{i}\omega,
  \\
  \lambda_ {6}^{\prime}=-\mathrm{i}\sqrt{ 3^{-\frac{\alpha}{2}}\alpha}=-\mathrm{i}\omega,
  \\
  \lambda_ {7}^{\prime}=\mathrm{i}\sqrt{ 3^{-\frac{\alpha}{2}}\alpha}=\mathrm{i}\omega,
  \\
  \lambda_ {8}^{\prime}=\mathrm{i}\sqrt{ 3^{-\frac{\alpha}{2}}\alpha}=\mathrm{i}\omega;
   \end{cases}
  \end{equation}
  
Then, for matrix $B_{3}$, we get 
  \begin{equation}
    B_{3}=\left(\begin{array}{cccc}
      0 & 0 & 1 & 0 \\
      0 & 0 & 0 & 1 \\
      3^{-\frac{\alpha}{2}}\alpha+\frac{1}{2} \times 3^{1-\frac{\alpha}{2}} \alpha^{2} & 0 & 0 & 2 \sqrt{3^{-\frac{\alpha}{2}}\alpha} \\
      0 & 3^{-\frac{\alpha}{2}}\alpha+\frac{1}{2} \times 3^{1-\frac{\alpha}{2}} \alpha^{2} & -2 \sqrt{3^{-\frac{\alpha}{2}}\alpha} & 0
      \end{array}\right).
  \end{equation}
  
  By calculating eigenvalues of $B_{3}$, we get 
  \begin{equation}
    \begin{cases}
  \lambda_ {9}^{\prime} =-\frac{\omega\sqrt{ (3\alpha -2\mathrm{i}\sqrt{6\alpha}-2 )}}{\sqrt{2}} ,
  \\
  \lambda_ {10}^{\prime}=\frac{\omega\sqrt{ (3\alpha -2\mathrm{i}\sqrt{6\alpha}-2 )}}{\sqrt{2}},
  \\
  \lambda_ {11}^{\prime}=-\frac{\omega\sqrt{ (3\alpha +2\mathrm{i}\sqrt{6\alpha}-2 )}}{\sqrt{2}},
  \\
  \lambda_ {12}^{\prime}=\frac{\omega\sqrt{ (3\alpha +2\mathrm{i}\sqrt{6\alpha}-2 )}}{\sqrt{2}}.
   \end{cases}
  \end{equation}
  
By analysing the above 12 eigenvalues, we find that for  $\alpha>0$, $\lambda_ {1}^{\prime},\lambda_ {2}^{\prime}$ are equal to zero; 
  
 When $\alpha\in(0,2)$, $\lambda_ {3}^{\prime},\lambda_ {4}^{\prime}$ are pure imaginary characteristic roots;
 
 When $\alpha=2$, $\lambda_ {3}^{\prime},\lambda_ {4}^{\prime}$ are zero roots;
 
 When $\alpha>2$, $\lambda_ {3}^{\prime},\lambda_ {4}^{\prime}$ are real characteristic roots;
  
 For  $\alpha>0$, $\lambda_ {5}^{\prime},\lambda_ {6}^{\prime},\lambda_ {7}^{\prime},\lambda_ {8}^{\prime}$
  are pure imaginary characteristic roots;

For  $\alpha>0$, $\lambda_ {9}^{\prime},\lambda_ {10}^{\prime},\lambda_ {11}^{\prime},\lambda_ {12}^{\prime}$ have nonzero real parts, since there are positive real parts, therefore the corresponding solutions are unstable.
  
 In conclusion, the relative equilibria of equal masses  of the equilateral triangle configuration under $\alpha$-homogeneous potential is unstable.
  
  \subsection{The stability of relative equilibria of four-body problem with equal mass}\label{sec: the stab of square homo}
  
  Let $ z \in \mathbb{R}^{8}$ be $z=\left(x_{1}, y_{1}, x_{2}, y_{2}, x_{3}, y_{3},x_{4}, y_{4}\right)^{T}$,
  we choose the  vertices of the square $$z_{0}=(1,0,0,1,-1,0,0,-1)^{T}$$
  as a critical point.
  
  Then we assume the corresponding relative equilibria are
  $$
    \mathbf{q}_{i}=\left(\begin{array}{c}
    \cos \left(\omega t+\frac{i-1}{2}\pi\right) \\
    \sin \left(\omega t+\frac{i+1}{2}\pi\right)
    \end{array}\right) , i=1,2,3,4.
  $$

  In $\alpha$-homogeneous four-body system with equal masses, the degree of $U$ and $I$ are  $\alpha$ and 2 respectively, we apply Euler formula
  \begin{equation}
    \omega^2=\frac{\alpha U(z_{0})}{2I(z_{0})}=(2^{-1-\alpha }+2^{-\frac{\alpha }{2}})\alpha .
  \end{equation}
  
  Then we obtain relative equilibria with initial configuration  $$z_{0}=(1,0,0,1,-1,0,0,-1)^{T}.$$
  
  We write the Hessian matrix $H_{2}:=D^2U(z_{0})$ as follows.
  $$
  H_{2}=
  \begin{pmatrix}
    A^{\prime}&C^{\prime}\\
    B^{\prime}&D^{\prime}
  \end{pmatrix}.
  $$
  where the matrices $A^{\prime},B^{\prime},C^{\prime},D^{\prime}$ are 
  $$
  A^{\prime}=
  \small{
  \begin{smallmatrix}
  \begin{pmatrix}
    \frac{2^{-\alpha} \alpha\left(1+\alpha+2^{1+\frac{\alpha}{2}} \alpha\right)}{4}& 0 & -2^{-2-\frac{\alpha}{2}} \alpha^{2} & 2^{-2-\frac{\alpha}{2}} \alpha(2+\alpha)\\
    0 & \frac{2^{-\alpha} \alpha\left(-1+2^{1+\frac{\alpha}{2}} \alpha\right)}{4} & 2^{-2-\frac{\alpha}{2}} \alpha(2+\alpha) & -2^{-2-\frac{\alpha}{2}} \alpha^{2}\\
    -2^{-2-\frac{\alpha}{2}} \alpha^{2} & 2^{-2-\frac{\alpha}{2}} \alpha(2+\alpha) & \frac{2^{-\alpha} \alpha\left(-1+2^{1+\frac{\alpha}{2}} \alpha\right)}{4} & 0\\
    2^{-2-\frac{\alpha}{2} \alpha(2+\alpha)} & -2^{-2-\frac{\alpha}{2}} \alpha^{2} & 0 & \frac{2^{-\alpha} \alpha\left(1+\alpha+2^{1+\frac{\alpha}{2}} \alpha\right)}{4}\\
  \end{pmatrix}
  \end{smallmatrix}
  },$$
  
  $$
  B^{\prime}=\begin{pmatrix}
    -2^{-2-\alpha} \alpha(1+\alpha) & 0 & -2^{-2-\frac{\alpha}{2}} \alpha^{2} & -2^{-2-\frac{\alpha}{2} \alpha(2+\alpha)}\\
    0 & 2^{-2-\alpha} \alpha & -2^{-2-\frac{\alpha}{2}} \alpha(2+\alpha) & -2^{-2-\frac{\alpha}{2}} \alpha^{2}\\
    -2^{-2-\frac{\alpha}{2}} \alpha^{2} & -2^{-2-\frac{\alpha}{2}} \alpha(2+\alpha) & 2^{-2-\alpha} \alpha & 0 \\
    -2^{-2-\frac{\alpha}{2}} \alpha(2+\alpha) & -2^{-2-\frac{\alpha}{2}} \alpha^{2} & 0 & -2^{-2-\alpha} \alpha(1+\alpha)
  \end{pmatrix},
  $$
  
  $$
  C^{\prime}=
  \begin{pmatrix}
    -2^{-2-\alpha} \alpha(1+\alpha) & 0 & -2^{-2-\frac{\alpha}{2}} \alpha^{2} & -2^{-2-\frac{\alpha}{2}} \alpha(2+\alpha) \\
    0 & 2^{-2-\alpha} \alpha & -2^{-2-\frac{\alpha}{2}} \alpha(2+\alpha) & -2^{-2-\frac{\alpha}{2}} \alpha^{2} \\
    -2^{-2-\frac{\alpha}{2}} \alpha^{2} & -2^{-2-\frac{\alpha}{2}} \alpha(2+\alpha) & 2^{-2-\alpha} \alpha & 0 \\
    -2^{-2-\frac{\alpha}{2} \alpha(2+\alpha)} & -2^{-2-\frac{\alpha}{2}} \alpha^{2} & 0 & -2^{-2-\alpha} \alpha(1+\alpha) \\
  \end{pmatrix},
  $$
  
  $$
  D^{\prime}=
  \small{
  \begin{smallmatrix}
  \begin{pmatrix}
    \frac{2^{-\alpha} \alpha\left(1+\alpha+2^{1+\frac{\alpha}{2}} \alpha\right)}{4} & 0 & -2^{-2-\frac{\alpha}{2}} \alpha^{2} & 2^{-2-\frac{\alpha}{2}} \alpha(2+\alpha) \\
    0 & \frac{2^{-\alpha} \alpha\left(-1+2^{1+\frac{\alpha}{2}} \alpha\right)}{4} & 2^{-2-\frac{\alpha}{2}} \alpha(2+\alpha) & -2^{-2-\frac{\alpha}{2}} \alpha^{2} \\
    -2^{-2-\frac{\alpha}{2}} \alpha^{2} & 2^{-2-\frac{\alpha}{2}} \alpha(2+\alpha) & \frac{2^{-\alpha} \alpha\left(-1+2^{1+\frac{\alpha}{2}} \alpha\right)}{4}& 0 \\
    2^{-2-\frac{\alpha}{2}} \alpha(2+\alpha) & -2^{-2-\frac{\alpha}{2}} \alpha^{2} & 0 & \frac{2^{-\alpha} \alpha\left(1+\alpha+2^{1+\frac{\alpha}{2}} \alpha\right)}{4}
  \end{pmatrix}
  \end{smallmatrix}.
  }$$
  
  We find that $\forall A\in H_{2}$, the equation $AH_{2}=H_{2}A$ is established with the effect of linear transformation,  we obtain the results.
  
  The trace of $H_{2}\mathscr{D}(A)$ is
  \begin{equation}
    \begin{aligned}
  \operatorname{Tr}\left(H_{2} \mathscr{D}(A)\right)&= 
  \lambda_{1} \chi_{1}(A)+\lambda_{2} \chi_{2}(A)+\lambda_{3} \chi_{3}(A)+\lambda_{4} \chi_{4}(A)+\lambda_{5} \chi_{5}(A)+\lambda_{6} \chi_{5}(A).
    \end{aligned}
  \end{equation}
  
Substituting the elements of $A$ to calculate the corresponding traces.
  \begin{equation}\label{equ:trace H_{4}}
  \begin{cases}
  \begin{aligned}
    \operatorname{Tr}\left(H_{2} \mathscr{D}(e)\right) &=\lambda_{1}+\lambda_{2}+\lambda_{3}+\lambda_{4}+2\left(\lambda_{5}+\lambda_{6}\right)=2^{-\alpha}\left(1+2^{2+\frac{\alpha}{2}}\right) \alpha^{2},
    \\
    \operatorname{Tr}\left(H_{2} \mathscr{D}(a)\right) &=\lambda_{1}+\lambda_{2}-\lambda_{3}-\lambda_{4} =0,
    \\
    \operatorname{Tr}\left(H_{2} \mathscr{D}\left(a^{2}\right)\right) &=\lambda_{1}+\lambda_{2}+\lambda_{3}+\lambda_{4}-2\left(\lambda_{5}+\lambda_{6}\right) =2^{-\alpha} \alpha^{2},
    \\
    \operatorname{Tr}\left(H_{2} \mathscr{D}(r)\right) &=\lambda_{1}-\lambda_{2}+\lambda_{3}-\lambda_{4}=2^{-\alpha} \alpha(2+\alpha),
     \\
    \operatorname{Tr}\left(H_{2} \mathscr{D}(a r)\right) &=\lambda_{1}-\lambda_{2}-\lambda_{3}+\lambda_{4}=2^{1-\frac{\alpha}{2}} \alpha(2+\alpha).
    \end{aligned}
  \end{cases}
  \end{equation}
  
To solve the eigenvalues, we combine the aforementioned equations.  
  \begin{equation}
  \begin{cases}
  \begin{aligned}
  \lambda_{1}&=2^{-1-\alpha}\left(1+2^{1+\frac{\alpha}{2}}\right) \alpha(1+\alpha),
  \\
  \lambda_{2}&=-2^{-1-\alpha}\left(1+2^{1+\frac{\alpha}{2}}\right) \alpha,
  \\
  \lambda_{3}&=-2^{-1-\alpha}\left(-1+2^{1+\frac{\alpha}{2}}-\alpha\right) \alpha, 
  \\
  \lambda_{4}&=2^{-1-\alpha} \alpha\left(-1+2^{1+\frac{\alpha}{2}}+2^{1+\frac{\alpha}{2}} \alpha\right),
  \\ 
  \lambda_{5}+\lambda_{6}&=2^{-\frac{\alpha}{2}}\alpha^{2}.
  \end{aligned}
  \end{cases}
  \end{equation}
  
  We calculate the characteristic polynomial $det(H_{2}+I)$ as follows. 
  \begin{equation}
  \begin{aligned}
   det(H_{2}+I)&=
  (\lambda_ {3}+1)( \lambda_ {4}+1)(\lambda_ {5}+1)^2( \lambda_ {6}+1)^2\\
  &=2^{4-5 \alpha}\left(2^{\alpha / 2}+\alpha^{2}\right)^{2} 
  \left [  16^{1+\alpha}+2^{3+2 \alpha}\left(-1-2^{1+\alpha}+2^{2+\frac{3}{2}}\right) \alpha^{2} \right ]\\
  &+2^{4-5 \alpha}\left(2^{\alpha / 2}+\alpha^{2}\right)^{2}\left [ 
  \left(1-3 \times 2^{2+\alpha}-2^{3+\frac{3 \alpha}{2}}-2^{3+\frac{5 \alpha}{2}}+
  2^{4+3 \alpha}+4^{1+\alpha}\right) \alpha^{4}\right ]\\
  &-2^{4-5 \alpha}\left(2^{\alpha / 2}+\alpha^{2}\right)^{2}\left [ 2\left(-1+5 \times 2^{1+\alpha}+
  2^{2+\frac{3 \alpha}{2}}-2^{3+2 \alpha}+2^{44 \frac{5 \alpha}{2}}\right) \alpha^{5}\right ]\\
  &-2^{4-5 \alpha}\left(2^{\alpha / 2}+\alpha^{2}\right)^{2}
  \left [ \left(1-2^{1+\frac{\alpha}{2}}-2^{4+\alpha}-2^{2+\frac{3 \alpha}{2}}+
  3 \times 2^{3+2 \alpha}\right) \alpha^{6}\right ]\\
  &-2^{4-5 \alpha}\left(2^{\alpha / 2}+\alpha^{2}\right)^{2}
  \left [ 2^{1+\frac{\alpha}{2}}\left(1+2^{1+\frac{\alpha}{2}}\right)^{2} \alpha^{7} +2^{3+2 \alpha}\left(1+2^{2+\alpha}\right) \alpha^{3}\right ].
  \end{aligned}
  \end{equation}

  Since there are double eigenvalues equal to 0, thus we assume $\lambda_{5}=0$, by combining the above equations, we solve   
  $$
  \lambda_{5}=0,\quad\lambda_{6}=2^{-\frac{\alpha}{2}}\alpha^{2}.
  $$
  
  In conclusion, all the eigenvalues of $H_{2}$ are 
  \begin{equation}
  \begin{cases}
  \lambda_{1}=2^{-1-\alpha}\left(1+2^{1+\frac{\alpha}{2}}\right) \alpha(1+\alpha),
  \\
  \lambda_{2}=-2^{-1-\alpha}\left(1+2^{1+\frac{\alpha}{2}}\right) \alpha, 
  \\
  \lambda_{3}=-2^{-1-\alpha}\left(-1+2^{1+\frac{\alpha}{2}}-\alpha\right) \alpha, 
  \\ 
  \lambda_{4}=2^{-1-\alpha} \alpha\left(-1+2^{1+\frac{\alpha}{2}}+2^{1+\frac{\alpha}{2}} \alpha\right), 
  \\
  \lambda_{5}=0,
  \\
  \lambda_{6}=2^{-\frac{\alpha}{2}}\alpha^{2}.
  \end{cases}
  \end{equation}
  Where the eigenvalues $\lambda_{5}$ and $\lambda_{6}$ are double characteristic roots, others are simple roots.
  
  By using our conclusion in section \ref*{sec: the stab of homo}, 

  For $\lambda_{1}$, we choose the eigenvector $u_{1}=\left\{1 , 0 , 0 , 1 ,-1 , 0 ,0 , -1\right\}^{T} $;

  For $\lambda_{2}$, we choose the eigenvector $u_{2}=\left\{0, 1 , -1 , 0 , 0 ,-1 , 1 ,0 \right\}^{T}$; 

  For $\lambda_{3}$, we choose the eigenvector $u_{3}=\left\{1 , 0 , 0 , -1 ,-1 , 0 ,0 , 1\right\}^{T} $;

  For $\lambda_{4}$, we choose the eigenvector $u_{4}=\left\{0, 1 , 1 , 0 , 0 ,-1 , -1 ,0 \right\}^{T}$; 
  
  Since $\lambda_{5}$ is a double eigenvalue, 
  thus we select two linearly independent eigenvectors,
  $u_{5}=\left\{1 ,0 ,1 ,0 ,1 ,0 ,1, 0\right\}^{T},
  u_{6}=\left\{0, 1, 0, 1, 0, 1, 0, 1 \right\}^{T}$; 

  The eigenvalue $\lambda_{6}$ is also a double eigenvalue, 
  thus we choose two linearly independent eigenvectors, 
  $u_{7}=\left\{1 ,0 ,-1 ,0 ,1 ,0 ,-1, 0\right\}^{T}
  $, $u_{8}=\left\{0, 1 ,0 ,-1 ,0 ,1 ,0 ,-1\right\}^{T}.
  $
  
Hence we choose these eigenvectors as a basis satisfying
    \begin{equation}
      \begin{split}
    J\begin{pmatrix}
      u_{1} &u_{2}
     \end{pmatrix}&=
     \left(\begin{array}{cc}
      1 & 0 \\
      0 & 1 \\
      0& -1 \\
      1& 0 \\
      -1& 0 \\
      0 & -1\\
      0& 1\\
      -1&0
      \end{array}\right)
      \begin{pmatrix}
      0 &1 &0 &0 &0 &0&0 &0\\
      -1&0 &0 &0 &0 &0&0 &0\\
       0&0 &0 &1 &0 &0&0 &0 \\
       0&0 &-1&0 &0 &0&0 &0 \\
       0&0 &0 &0 &0 &1&0 &0\\
       0&0 &0 &0 &-1&0&0 &0\\
       0&0 &0 &0 &0 &0&0 &1\\
       0&0 &0 &0 &0 &0&-1 &0
     \end{pmatrix}
     =
     \left(\begin{array}{cc}
     0 & 1 \\
     -1 & 0 \\
     1 & 0 \\
     0 & 1 \\
     0 & -1\\
     1 & 0\\
     -1& 0\\
     0&-1
     \end{array}\right) \\
    & =
    \begin{pmatrix}
      0 &1 &0 &0 &0 &0&0 &0\\
      -1&0 &0 &0 &0 &0&0 &0\\
       0&0 &0 &1 &0 &0&0 &0 \\
       0&0 &-1&0 &0 &0&0 &0 \\
       0&0 &0 &0 &0 &1&0 &0\\
       0&0 &0 &0 &-1&0&0 &0\\
       0&0 &0 &0 &0 &0&0 &1\\
       0&0 &0 &0 &0 &0&-1 &0
     \end{pmatrix}
     \left(\begin{array}{cc}
     1 & 0 \\
     0 & 1 \\
     0& -1 \\
     1& 0 \\
     -1& 0 \\
     0 & -1\\
     0& 1\\
     -1&0
     \end{array}\right)=\begin{pmatrix}
      u_{1} &u_{2}
     \end{pmatrix}J.
    \end{split}
    \end{equation}
  
   Similarly, we find for $\lambda_{3}$ and $\lambda_{4}$, we have
    $
      J \begin{pmatrix}
        u_{3} &u_{4}
       \end{pmatrix}=\begin{pmatrix}
        u_{3} &u_{4}
       \end{pmatrix}J
    $, for $\lambda_{5}=0$, we have
    $
      J \begin{pmatrix}
        u_{5} &u_{6}
       \end{pmatrix}=\begin{pmatrix}
        u_{5} &u_{6}
       \end{pmatrix}J
    $, for $\lambda_{6}=2^{-\frac{\alpha}{2}}\alpha^{2}$,
    $
      J \begin{pmatrix}
        u_{7} &u_{8}
       \end{pmatrix}=\begin{pmatrix}
        u_{7} &u_{8}
       \end{pmatrix}J
  $.
  
  In the new coordinate $E_{\lambda}=\left \{ V_{i}, V_{j}\right \}$, let $y_{i}=(x_{i},\dot{x_{i}})$, then the linearized second order equation is 
  $\dot{y_{i}}=B_{i}y_{i},\,i=1,2,3,4$,
  where $B_{i}$ satisfy
  
  \begin{equation}
  B_{i}=\left(\begin{array}{cc}
    0 & I_{2} \\
    \omega^{2} I_{2}+\frac{1}{m_{i}} \Lambda_{i} & 2 \omega J
    \end{array}\right),i=1,2,3,4,
  \end{equation}
  the $\Lambda_{i}$, i=1,2,3,4 are 
  $$
  \begin{array}{c}
    \Lambda_{1}=\left(\begin{array}{cc}
    \lambda_{1} & 0 \\
    0 & \lambda_{2}
    \end{array}\right), \quad \Lambda_{2}=\left(\begin{array}{cc}
    \lambda_{3} & 0 \\
    0 & \lambda_{4}
    \end{array}\right), \quad 
    \Lambda_{3}=\left(\begin{array}{cc}
    \lambda_{5} & 0 \\
    0 & \lambda_{5}
    \end{array}\right),\quad 
    \Lambda_{4}=\left(\begin{array}{cc}
    \lambda_{6} & 0 \\
    0 & \lambda_{6}
    \end{array}\right)
  \end{array}.
  $$
  
Thus we calculate all the eigenvalues of equations $\dot{y_{i}}=B_{i}y_{i},\,i=1,2,3,4$.
  
  Firstly, we calculate the matrix $B_{1}$.
  \begin{equation}
    \begin{smallmatrix}
    \begin{aligned}
  B_{1} &=
  \left(\begin{array}{cccc}
    0 & 0 & 1 & 0 \\
    0 & 0 & 0 & 1 \\
    \omega^{2}+\lambda_{1} & 0 & 0 & 2 \omega \\
    0 & \omega^{2}+\lambda_{2} & -2 \omega & 0
    \end{array}\right) ,
  \end{aligned}
  \end{smallmatrix}
  \end{equation}
  where $\omega^2=(2^{-1-\alpha }+2^{-\frac{\alpha }{2}})\alpha$, 
  $\lambda_{1}=2^{-1-\alpha}\left(1+2^{1+\frac{\alpha}{2}}\right) \alpha(1+\alpha)$, $\lambda_{2}=-2^{-1-\alpha}\left(1+2^{1+\frac{\alpha}{2}}\right) \alpha.$
  
  By calculating eigenvalues of $B_{1}$, we get 
  
  \begin{equation}
    \begin{cases}
  \lambda_ {1}^{\prime} =0,
  \\
  \lambda_ {2}^{\prime}=0,\\
  \lambda_ {3}^{\prime}=-\mathrm{i} \sqrt{ (2^{-1-\alpha }+2^{-\frac{\alpha }{2}})\alpha(2-\alpha)}=-\mathrm{i}\omega\sqrt{(2-\alpha)},
  \\
  \lambda_ {4}^{\prime}=\mathrm{i} \sqrt{ (2^{-1-\alpha }+2^{-\frac{\alpha }{2}})\alpha(2-\alpha)}=\mathrm{i}\omega\sqrt{(2-\alpha)}.
   \end{cases}
  \end{equation}
  
  Then, we calculate the matrix $B_{2}$.
  \begin{equation}
  B_{2}=
  \left(\begin{array}{cccc}
    0 & 0 & 1 & 0 \\
    0 & 0 & 0 & 1 \\
    \omega^2+\lambda_{3} & 0 & 0 & 2 \omega \\
    0 & \omega^2+\lambda_{4}   & -2 \omega & 0
    \end{array}\right),
  \end{equation}
  where $\lambda_{3}=-2^{-1-\alpha}\left(-1+2^{1+\frac{\alpha}{2}}-\alpha\right) \alpha$, $\lambda_{4}=2^{-1-\alpha} 
  \alpha\left(-1+2^{1+\frac{\alpha}{2}}+2^{1+\frac{\alpha}{2}} \alpha\right).$
  
  By calculating eigenvalues of $B_{2}$, we get 
  \begin{equation}
    \begin{cases}
  \lambda_ {5}^{\prime} =\lambda_ {6}^{\prime}= -\sqrt{-\frac{2-\alpha }{2}\omega ^2 +2^{-2-\alpha} \alpha \sqrt{m}},
  \\
  \lambda_ {7}^{\prime}= \lambda_ {8}^{\prime}=-\sqrt{(-\frac{2-\alpha }{2}\omega ^2 -2^{-2-\alpha} \alpha \sqrt{m}},
   \end{cases}
  \end{equation}
  where $m=-2^{4+\frac{\alpha}{2}}(1+3\alpha )+2^{4+\alpha}(1-\alpha )+(\alpha-2)^{2}-2^{2+\frac{\alpha}{2}} \alpha^{2}(1-2^{\frac{\alpha }{2} }).$
 
  We calculate the matrix $B_{3}$.
  \begin{equation}
    B_{3}=
      \left(\begin{array}{cccc}
      0 & 0 & 1 & 0 \\
      0 & 0 & 0 & 1 \\
      (2^{-1-\alpha }+2^{-\frac{\alpha }{2}})\alpha+4^{-\alpha} & 0 & 0 & 2\omega \\
      0 & (2^{-1-\alpha }+2^{-\frac{\alpha }{2}})\alpha+4^{-\alpha} & -2 \omega & 0
      \end{array}\right).
  \end{equation}
  
Calculating the eigenvalues of $B_{3}$, we get 
  \begin{equation}
    \begin{cases}
  \lambda_ {9}^{\prime} =-\mathrm{i}\sqrt{(2^{-1-\alpha }+2^{-\frac{\alpha }{2}})\alpha}=-\mathrm{i}\omega,
  \\
  \lambda_ {10}^{\prime}=-\mathrm{i}\sqrt{(2^{-1-\alpha }+2^{-\frac{\alpha }{2}})\alpha}=-\mathrm{i}\omega,
  \\
  \lambda_ {11}^{\prime}=\mathrm{i}\sqrt{(2^{-1-\alpha }+2^{-\frac{\alpha }{2}})\alpha}=\mathrm{i}\omega,
  \\
  \lambda_ {12}^{\prime}=\mathrm{i}\sqrt{(2^{-1-\alpha }+2^{-\frac{\alpha }{2}})\alpha}=\mathrm{i}\omega.
   \end{cases}
  \end{equation}
  
  We calculate the matrix $B_{4}$ to obtain
  \begin{equation}
    B_{4}=\left(\begin{array}{cccc}
      0 & 0 & 1 & 0 \\
      0 & 0 & 0 & 1 \\
      (2^{-1-\alpha }+2^{-\frac{\alpha }{2}})\alpha+2^{-\frac{\alpha}{2}}\alpha^{2} & 0 & 0 & 2 \omega \\
      0 & (2^{-1-\alpha }+2^{-\frac{\alpha }{2}})\alpha+2^{-\frac{\alpha}{2}}\alpha^{2} & -2 \omega & 0
      \end{array}\right).
  \end{equation}
  
Calculating eigenvalues of $B_{4}$, we get 
  \begin{equation}
    \begin{cases}
  \lambda_ {13}^{\prime} =-2^{-\frac{1}{2}-\alpha} \sqrt{-2^{\alpha} \alpha-2^{1+\frac{3 \alpha}{2}} \alpha+2^{1+\frac{3 \alpha}{2}} \alpha^{2}-2 \sqrt{2} \sqrt{-2^{\frac{5\alpha }{2}}\left(1+2^{1+\frac{\alpha}{2}}\right) \alpha^{3}}},
  \\
  \lambda_ {14}^{\prime}=2^{-\frac{1}{2}-\alpha} \sqrt{-2^{\alpha} \alpha-2^{1+\frac{3 \alpha}{2}} \alpha+2^{1+\frac{3 \alpha}{2}} \alpha^{2}-2 \sqrt{2} \sqrt{-2^{\frac{5\alpha }{2}}\left(1+2^{1+\frac{\alpha}{2}}\right) \alpha^{3}}},
  \\
  \lambda_ {15}^{\prime}=-2^{-\frac{1}{2}-\alpha} \sqrt{-2^{\alpha} \alpha-2^{1+\frac{3 \alpha}{2}} \alpha+2^{1+\frac{3 \alpha}{2}} \alpha^{2}+2 \sqrt{2} \sqrt{-2^{\frac{5\alpha }{2}}\left(1+2^{1+\frac{\alpha}{2}}\right) \alpha^{3}}},
  \\
  \lambda_ {16}^{\prime}=2^{-\frac{1}{2}-\alpha} \sqrt{-2^{\alpha} \alpha-2^{1+\frac{3 \alpha}{2}} \alpha+2^{1+\frac{3 \alpha}{2}} \alpha^{2}+2 \sqrt{2} \sqrt{-2^{\frac{5\alpha }{2}}\left(1+2^{1+\frac{\alpha}{2}}\right) \alpha^{3}}}.
   \end{cases}
  \end{equation}

By analysing the above 16 eigenvalues, we find that for  $\alpha>0$,
eigenvalues $\lambda_ {1}^{\prime}, \lambda_ {2}^{\prime}$ are trivial solutions;

When $\alpha\in(0,2)$, $\lambda_ {3}^{\prime},\lambda_ {4}^{\prime}$ are pure imaginary characteristic roots;
 
When $\alpha=2$,$\lambda_ {3}^{\prime},\lambda_ {4}^{\prime}$ are zero roots;
 
When $\alpha>2$, $\lambda_ {3}^{\prime},\lambda_ {4}^{\prime}$ are real characteristic roots;

For  $\alpha>0$, eigenvalues $ \lambda_ {9}^{\prime}, \lambda_ {10}^{\prime}, \lambda_ {11}^{\prime}, \lambda_ {12}^{\prime}$ 
are pure imaginary roots;

For  $\alpha>0$, eigenvalues $\lambda_ {5}^{\prime}, \lambda_ {6}^{\prime}, \lambda_ {7}^{\prime}, \lambda_ {8}^{\prime}, \lambda_ {13}^{\prime}, \lambda_ {14}^{\prime}, \lambda_ {15}^{\prime}, \lambda_ {16}^{\prime}$ 
all have nonzero real parts.

 In conclusion, the relative equilibria of equal mass of the square configuration under $\alpha$-homogeneous potential is unstable.

  
  \section{The linear stability of regular polygon under quasi-homogeneous potential with equal mass}\label{sec: the stab of quasihomo}

  \label{sec:examples}
  \subsection{The linear stability of equilateral triangle problem under the Manev potential with equal mass}\label{sec: the stab of tria Manev}

  In this section, we consider the stability of equal masses of equilateral triangle problem in $\mathbb{R}^2$, $U$ satisfies
      \begin{equation}
          U=\sum\limits_{1 \leq i<j \leq 3} \frac{1}{\left [  \left(x_{i}-x_{j}\right)^{2}+\left(y_{i}-y_{j}\right)^{2}\right ]^{\frac{1 }{2} } }+\frac{1}{ \left(x_{i}-x_{j}\right)^{2}+\left(y_{i}-y_{j}\right)^{2}}.  
      \end{equation}
  
  Let the center of mass be at the origin, and
  $ z \in \mathbb{R}^{6}$ be $z=\left(x_{1}, y_{1}, x_{2}, y_{2}, x_{3}, y_{3}\right)^{T}$,
  we choose the vertices of equilateral triangle $$z_{0}=(1,0,-\frac{1}{2}, \frac{\sqrt{3}}{2},-\frac{1}{2},-\frac{\sqrt{3}}{2})^{T},$$ to solve the equation $\nabla (U+\omega ^2I)=0, \omega \ne 0$,
  we obtain $\omega^2=\frac{2}{3}+\frac{1}{\sqrt{3}}$.
  
  At the point $z_{0}$, the Hessian matrix of $H_{3}$ at $z_{0}$ is  
  $$
  H_{3}=
  \begin{pmatrix}
      \begin{smallmatrix}
          \frac{1}{9}(16+5 \sqrt{3}) & 0 & -\frac{8}{9}-\frac{5}{6 \sqrt{3}} & \frac{1}{18}(9+8 \sqrt{3}) & -\frac{8}{9}-\frac{5}{6 \sqrt{3}} & \frac{1}{18}(-9-8 \sqrt{3}) \\
          0 & -\frac{1}{3 \sqrt{3}} & \frac{1}{18}(9+8 \sqrt{3}) & \frac{1}{6 \sqrt{3}} & \frac{1}{18}(-9-8 \sqrt{3}) & \frac{1}{6 \sqrt{3}} \\
          -\frac{8}{9}-\frac{5}{6 \sqrt{3}} & \frac{1}{18}(9+8 \sqrt{3}) & \frac{1}{18}(8+\sqrt{3}) & \frac{1}{18}(-9-8 \sqrt{3}) & \frac{2}{9}(2+\sqrt{3}) & 0 \\
          \frac{1}{18}(9+8 \sqrt{3}) & \frac{1}{6 \sqrt{3}} & \frac{1}{18}(-9-8 \sqrt{3}) & \frac{4}{3}+\frac{7}{6 \sqrt{3}} & 0 & -\frac{4}{9}(3+\sqrt{3}) \\
          -\frac{8}{9}-\frac{5}{6 \sqrt{3}} & \frac{1}{18}(-9-8 \sqrt{3}) & \frac{2}{9}(2+\sqrt{3}) & 0 & \frac{1}{18}(8+\sqrt{3}) & \frac{1}{18}(9+8 \sqrt{3}) \\
          \frac{1}{18}(-9-8 \sqrt{3}) & \frac{1}{6 \sqrt{3}} & 0 & -\frac{4}{9}(3+\sqrt{3}) & \frac{1}{18}(9+8 \sqrt{3}) & \frac{4}{3}+\frac{7}{6 \sqrt{3}}        
  \end{smallmatrix}    
  \end{pmatrix}.
  $$
  
Applying the previous theory, we compute the traces of $H_{3}\mathscr{D}(A)$ for different irreducible group representations, and we find that the eigenvalues of the matrix $H_{3}$ are as follows.
  
  \begin{equation}
   \begin{split}
     \begin{cases}
      \lambda_ {1}=\frac{4}{3}(3+\sqrt{3}),
      \\
      \lambda_ {2}=-\frac{2}{3}(2+\sqrt{3}),
      \\
      \lambda_ {3}=0,
      \\
      \lambda_ {4}=\frac{1}{3}(4+\sqrt{3}).
      \end{cases}
    \end{split}
  \end{equation}
  
  Then linearize the second order equation 
  \begin{equation}\label{equ:linearied second equation2}
    \ddot{x}_{i}=2 \omega J \dot{x_{i}}+\omega^{2} x_{i}+\nabla_{x_{i}} U, x\in \mathbb{R}^2.
  \end{equation} 

  In the new coordinate $E_{\lambda}=\left \{ V_{i}, V_{j}\right \}$, let $y_{i}=(x_{i},\dot{x_{i}})$, therefore the linearized second order equations are 
  $\dot{y_{i}}=B_{i}y_{i},\,i=1,2,3$.
  where $B_{i},i=1,2,3,$ satisfy
  \begin{equation}
  B_{i}=\left(\begin{array}{cc}
    0 & I_{2} \\
    \omega^{2} I_{2}+\frac{1}{m_{i}} \Lambda_{i} & 2 \omega J
    \end{array}\right),i=1,2,3.
  \end{equation}
  The $\Lambda_{i}$, i=1,2,3 are 
  $$
  \begin{array}{c}
    \Lambda_{1}=\left(\begin{array}{cc}
    \lambda_{1} & 0 \\
    0 & \lambda_{2}
    \end{array}\right), \quad \Lambda_{2}=\left(\begin{array}{cc}
    \lambda_{3} & 0 \\
    0 & \lambda_{3}
    \end{array}\right), \quad 
    \Lambda_{3}=\left(\begin{array}{cc}
    \lambda_{4} & 0 \\
    0 & \lambda_{4}
    \end{array}\right)
  \end{array}.
  $$
  
  By calculating the eigenvalues of matrix $B_{1}$, we get 
  \begin{equation}
      \begin{cases}
    \lambda_ {1}^{\prime} =-\sqrt{\frac{1}{6}(5-\sqrt{3}+2 \sqrt{49+18 \sqrt{3}})},
    \\
    \lambda_ {2}^{\prime}=\sqrt{\frac{1}{6}(5-\sqrt{3}+2 \sqrt{49+18 \sqrt{3}})},
    \\
    \lambda_ {3}^{\prime}=-\mathrm{i}\sqrt{\frac{1}{6}(5-\sqrt{3}+2 \sqrt{49+18 \sqrt{3}})},
    \\
    \lambda_ {4}^{\prime}=\mathrm{i}\sqrt{\frac{1}{6}(5-\sqrt{3}+2 \sqrt{49+18 \sqrt{3}})};
      \end{cases}
    \end{equation}
    
  By calculating the eigenvalues of $B_{2}$, we get 
  \begin{equation}
      \begin{cases}
    \lambda_ {5}^{\prime} =\lambda_ {6}^{\prime}=-\mathrm{i} \sqrt{\frac{1}{2}(1+\sqrt{3})},
    \\
    \lambda_ {7}^{\prime}=\lambda_ {8}^{\prime}=\mathrm{i} \sqrt{\frac{1}{2}(1+\sqrt{3})};
     \end{cases}
    \end{equation}
  
  By calculating the eigenvalues of $B_{3}$, we get 
    \begin{equation}
      \begin{cases}
    \lambda_ {9}^{\prime} =-\sqrt{\frac{5}{6}-\frac{1}{2 \sqrt{3}}-\frac{1}{6} \mathrm{i} \sqrt{168+120 \sqrt{3}}},
    \\
    \lambda_ {10}^{\prime}=\sqrt{\frac{5}{6}-\frac{1}{2 \sqrt{3}}-\frac{1}{6} \mathrm{i} \sqrt{168+120 \sqrt{3}}},
    \\
    \lambda_ {11}^{\prime}=-\sqrt{\frac{5}{6}-\frac{1}{2 \sqrt{3}}+\frac{1}{6} \mathrm{i} \sqrt{168+120 \sqrt{3}}},
    \\
    \lambda_ {12}^{\prime}=\sqrt{\frac{5}{6}-\frac{1}{2 \sqrt{3}}+\frac{1}{6} \mathrm{i} \sqrt{168+120 \sqrt{3}}}.
     \end{cases}
    \end{equation}
  
  By analysing the above 12 eigenvalues, we conclude that 
  $\lambda_ {1}^{\prime},\lambda_ {2}^{\prime}$ are real roots; 
  $\lambda_ {3}^{\prime},\lambda_ {4}^{\prime}, \lambda_ {5}^{\prime},\lambda_ {6}^{\prime},\lambda_ {7}^{\prime},\lambda_ {8}^{\prime}$
  are pure imaginary characteristic roots, 
  $\lambda_ {9}^{\prime},\lambda_ {10}^{\prime},\lambda_ {11}^{\prime},\lambda_ {12}^{\prime}$ are imaginary roots with nonzero real parts, thus the corresponding solutions are unstable.
  
  In conclusion, the relative equilibria of the  equilateral triangle configuration  with equal mass under the Manev potential is unstable.

  \subsection{The stability of equilateral triangle problem under the Schwarzschild potential with equal mass}\label{sec: the stab of tria Schw}

  In this section, we consider the stability with equal mass of equilateral triangle configuration under the Schwarzschild potential in $\mathbb{R}^2$, $U$ satisfies
      \begin{equation}
        U=\sum\limits_{1 \leq i<j \leq 3} \frac{1}{\left [  \left(x_{i}-x_{j}\right)^{2}+\left(y_{i}-y_{j}\right)^{2}\right ]^{\frac{3 }{2} } }+\frac{1}{\sqrt{ \left(x_{i}-x_{j}\right)^{2}+\left(y_{i}-y_{j}\right)^{2}}}.
      \end{equation}

      Let the center of mass be at the origin, and
      $ z \in \mathbb{R}^{6}$ be $z=\left(x_{1}, y_{1}, x_{2}, y_{2}, x_{3}, y_{3}\right)^{T}$,
      we choose the vertices of equilateral triangle $$z_{0}=(1,0,-\frac{1}{2}, \frac{\sqrt{3}}{2},-\frac{1}{2},-\frac{\sqrt{3}}{2})^{T}$$ as an initial point to solve the equation $\nabla (U+\omega ^2I)=0, \omega \ne 0$,
      we obtain $\omega^2=\frac{2}{\sqrt{3}}$.

  At the point $z_{0}$, the Hessian matrix of $H_{4}$ at $z_{0}$ is  
  $$
  H_{4}=
  \begin{pmatrix}
        \frac{8}{3 \sqrt{3}} & 0 & -\frac{4}{3 \sqrt{3}} & \frac{2}{3} & -\frac{4}{3 \sqrt{3}} & -\frac{2}{3} \\
        0 & 0 & \frac{2}{3} & 0 & -\frac{2}{3} & 0 \\
        -\frac{4}{3 \sqrt{3}} & \frac{2}{3} & \frac{2}{3 \sqrt{3}} & -\frac{2}{3} & \frac{2}{3 \sqrt{3}} & 0 \\
        \frac{2}{3} & 0 & -\frac{2}{3} & \frac{2}{\sqrt{3}} & 0 & -\frac{2}{\sqrt{3}} \\
        -\frac{4}{3 \sqrt{3}} & -\frac{2}{3} & \frac{2}{3 \sqrt{3}} & 0 & \frac{2}{3 \sqrt{3}} & \frac{2}{3} \\
        -\frac{2}{3} & 0 & 0 & -\frac{2}{\sqrt{3}} & \frac{2}{3} & \frac{2}{\sqrt{3}}          
  \end{pmatrix}.
  $$
  
  As before, we calculate the traces of $H_{4}\mathscr{D}(A)$ to solve the eigenvalues of matrix $H_{4}$
  \begin{equation}
    \begin{split}
      \begin{cases}
          \lambda_ {1}=\frac{29}{6\sqrt{3}},
          \\
          \lambda_ {2}=-\frac{19}{6\sqrt{3}},
          \\
          \lambda_ {3}=0,
          \\
          \lambda_ {4}=-\frac{1}{3\sqrt{3}}.
          \end{cases}
    \end{split}
  \end{equation}
  
  Then we linearize the second-order equation 
  \begin{equation}\label{equ:linearied second equation3}
    \ddot{x}_{i}=2 \omega J \dot{x_{i}}+\omega^{2} x_{i}+\nabla_{x_{i}} U, x\in \mathbb{R}^2.
  \end{equation} 
  In the new coordinate $E_{\lambda}=\left \{ V_{i}, V_{j}\right \}$, let $y_{i}=(x_{i},\dot{x_{i}})$, the linearized second order equations are 
  $\dot{y_{i}}=B_{i}y_{i},\,i=1,2,3$.

 The matrices $B_{i},i=1,2,3,$ satisfy
  \begin{equation}
  B_{i}=\left(\begin{array}{cc}
    0 & I_{2} \\
    \omega^{2} I_{2}+\frac{1}{m_{i}} \Lambda_{i} & 2 \omega J
    \end{array}\right),i=1,2,3.
  \end{equation}
  The $\Lambda_{i}$, i=1,2,3 are 
  $$
  \begin{array}{c}
    \Lambda_{1}=\left(\begin{array}{cc}
    \lambda_{1} & 0 \\
    0 & \lambda_{2}
    \end{array}\right), \quad \Lambda_{2}=\left(\begin{array}{cc}
    \lambda_{3} & 0 \\
    0 & \lambda_{3}
    \end{array}\right), \quad 
    \Lambda_{3}=\left(\begin{array}{cc}
    \lambda_{4} & 0 \\
    0 & \lambda_{4}
    \end{array}\right)
  \end{array}.
  $$

  By calculating eigenvalues of matrix $B_{1}$, we get 
  \begin{equation}
        \lambda_ {1}^{\prime} =\lambda_ {2}^{\prime}=\lambda_ {3}^{\prime}=\lambda_ {4}^{\prime}=0;
    \end{equation}
    
  By calculating eigenvalues of $B_{2}$, we get 
  \begin{equation}
      \begin{cases}
          \lambda_ {5}^{\prime} =\lambda_ {6}^{\prime}=-(\sqrt{2})\times (3^{-{\frac{1}{4}}})\mathrm{i},
          \\
          \lambda_ {7}^{\prime}=\lambda_ {8}^{\prime}=(\sqrt{2})\times (3^{-{\frac{1}{4}}})\mathrm{i};
     \end{cases}
    \end{equation}
  
  By calculating eigenvalues of $B_{3}$, we get 
    \begin{equation}
      \begin{cases}
          \lambda_ {9}^{\prime} =-2(-\frac{1}{3})^{\frac{1}{4}},
          \\
          \lambda_ {10}^{\prime}=2(-\frac{1}{3})^{\frac{1}{4}},
          \\
          \lambda_ {11}^{\prime}=-2(-\frac{1}{3})^{\frac{3}{4}},
          \\
          \lambda_ {12}^{\prime}=2(-\frac{1}{3})^{\frac{3}{4}}.
     \end{cases}
    \end{equation}
  
    By analysing the above 12 eigenvalues, we find that 
   $\lambda_ {1}^{\prime},\lambda_ {2}^{\prime},\lambda_ {3}^{\prime},\lambda_ {4}^{\prime}$ are zero roots,
  $ \lambda_ {5}^{\prime},\lambda_ {6}^{\prime},\lambda_ {7}^{\prime},\lambda_ {8}^{\prime}$ are pure imaginary characteristic roots,
  $ \lambda_ {9}^{\prime},\lambda_ {10}^{\prime},\lambda_ {11}^{\prime},\lambda_ {12}^{\prime}$ are imaginary characteristic roots with nonzero real parts, thus the corresponding solutions are unstable.
  
  In conclusion, the relative equilibria of the equilateral triangle configuration with equal mass under the Schwarzschild potential is unstable.
  

  \subsection{The stability of square problem under the Manev potential with equal mass}\label{sec: the stab of square Manev}

  In this section, we consider the stability of equal mass of equilateral triangle problem in $\mathbb{R}^2$, $U$ satisfies
      \begin{equation}
          U=\sum\limits_{1 \leq i<j \leq 4} \frac{1}{\left [  \left(x_{i}-x_{j}\right)^{2}+\left(y_{i}-y_{j}\right)^{2}\right ]^{\frac{1 }{2} } }+\frac{1}{ \left(x_{i}-x_{j}\right)^{2}+\left(y_{i}-y_{j}\right)^{2}}.  
      \end{equation}
  
  Let the center of mass be at the origin, and
  $z \in \mathbb{R}^{8}$ be $z=\left(x_{1}, y_{1}, x_{2}, y_{2}, x_{3}, y_{3}, x_{4}, y_{4}\right)^{T}$ 
  we choose the vertices of square $$z_{0}=(1,0,0, 1,-1,0,0,-1)^{T}$$ to solve the equation $\nabla (U+\omega ^2I)=0, \omega \ne 0$,
  we obtain $\omega^2=\frac{2}{3}+\frac{1}{\sqrt{3}}$.
  
  At the point $z_{0}$, the Hessian matrix of $H_{5}$ at $z_{0}$ is  
  $$
  H_{5}=
  \begin{pmatrix}
      \begin{smallmatrix}
          \frac{1}{8}(13+2 \sqrt{2}) & 0 & \frac{1}{8}(-4-\sqrt{2}) & 1+\frac{3}{4 \sqrt{2}} & -\frac{5}{8} & 0 & \frac{1}{8}(-4-\sqrt{2}) & -1-\frac{3}{4 \sqrt{2}}   \\
  0 & \frac{1}{4}(3+\sqrt{2}) & 1+\frac{3}{4 \sqrt{2}} & \frac{1}{8}(-4-\sqrt{2}) & 0 & \frac{1}{4} &  -1-\frac{3}{4 \sqrt{2}} & \frac{1}{8}(-4-\sqrt{2}) \\
  \frac{1}{8}(-4-\sqrt{2}) & 1+\frac{3}{4 \sqrt{2}} & \frac{1}{4}(3+\sqrt{2}) & 0 & \frac{1}{8}(-4-\sqrt{2}) & -1-\frac{3}{4 \sqrt{2}} &  \frac{1}{4} & 0 \\
  1+\frac{3}{4 \sqrt{2}} & \frac{1}{8}(-4-\sqrt{2}) & 0 & \frac{1}{8}(13+2 \sqrt{2}) & -1-\frac{3}{4 \sqrt{2}}& \frac{1}{8}(-4-\sqrt{2}) & 0&-\frac{5}{8}\\
  -\frac{5}{8} & 0 & \frac{1}{8}(-4-\sqrt{2}) & -1-\frac{3}{4 \sqrt{2}} & \frac{1}{8}(13+2 \sqrt{2}) & 0 & \frac{1}{8}(-4-\sqrt{2}) & 1+\frac{3}{4\sqrt{2}} \\
  0 & \frac{1}{4} & -1-\frac{3}{4 \sqrt{2}} & \frac{1}{8}(-4-\sqrt{2}) & 0 & \frac{1}{4}(3+\sqrt{2})&  1+\frac{3}{4 \sqrt{2}} & \frac{1}{8}(-4-\sqrt{2}) \\
  \frac{1}{8}(-4-\sqrt{2}) & -1-\frac{3}{4 \sqrt{2}} & \frac{1}{4} & 0 & \frac{1}{8}(-4-\sqrt{2}) & 1+\frac{3}{4 \sqrt{2}} & \frac{1}{4}(3+\sqrt{2}) & 0 \\
  -1-\frac{3}{4 \sqrt{2}} & \frac{1}{8}(-4-\sqrt{2}) & 0 & -\frac{5}{8} & 1+\frac{3}{4 \sqrt{2}} & \frac{1}{8}(-4-\sqrt{2}) & 0 & \frac{1}{8}(13+2 \sqrt{2}) 
  \end{smallmatrix}
  \end{pmatrix}.
  $$

To solve the eigenvalues of matrix $H_{5}$, we calculate the traces of $H_{5}\mathscr{D}(A)$ for different irreducible group representations.
  
  \begin{equation}
   \begin{split}
     \begin{cases}
      \lambda_ {1}=\frac{17}{4}+\sqrt{2},
      \\
      \lambda_ {2}=-\frac{3}{2}-\frac{1}{\sqrt{2}},
      \\
      \lambda_ {3}=\frac{1}{4}-\frac{1}{\sqrt{2}},
      \\
      \lambda_ {4}=\frac{5}{2}+\sqrt{2},
      \\
      \lambda_ {5}=0,
      \\
      \lambda_ {6}=2+\frac{1}{\sqrt{2}}.
      \end{cases}
    \end{split}
  \end{equation}
  
  Furthermore, we linearize the second-order equation 
  \begin{equation}\label{equ:linearied second equation4}
    \ddot{x}_{i}=2 \omega J \dot{x_{i}}+\omega^{2} x_{i}+\nabla_{x_{i}} U, x\in \mathbb{R}^2.
  \end{equation} 
  In the new coordinate $E_{\lambda}=\left \{ V_{i}, V_{j}\right \}$, let $y_{i}=(x_{i},\dot{x_{i}})$, the linearized second order equation are 
  $\dot{y_{i}}=B_{i}y_{i},\,i=1,2,3,4$.

  The $B_{i}, i=1,2,3,4,$ satisfy
  \begin{equation}
  B_{i}=\left(\begin{array}{cc}
    0 & I_{2} \\
    \omega^{2} I_{2}+\frac{1}{m_{i}} \Lambda_{i} & 2 \omega J
    \end{array}\right),i=1,2,3,4.
  \end{equation}
  The $\Lambda_{i}$, i=1,2,3,4, are 
  $$
  \begin{array}{c}
    \Lambda_{1}=\left(\begin{array}{cc}
    \lambda_{1} & 0 \\
    0 & \lambda_{2}
    \end{array}\right), \quad \Lambda_{2}=\left(\begin{array}{cc}
    \lambda_{3} & 0 \\
    0 & \lambda_{4}
    \end{array}\right), \quad 
    \Lambda_{3}=\left(\begin{array}{cc}
    \lambda_{5} & 0 \\
    0 & \lambda_{5}
    \end{array}\right), \quad 
    \Lambda_{4}=\left(\begin{array}{cc}
    \lambda_{6} & 0 \\
    0 & \lambda_{6}
    \end{array}\right)
  \end{array}.
  $$
  By calculating the eigenvalues of matrix $B_{1}$, we get 
  \begin{equation}
      \begin{cases}
          \lambda_ {1}^{\prime} =0,
          \\
          \lambda_ {2}^{\prime}=0,
          \\
          \lambda_ {3}^{\prime}=-\frac{1}{2} \mathrm{i} \sqrt{1+2 \sqrt{2}},
          \\
          \lambda_ {4}^{\prime}=\frac{1}{2} \mathrm{i} \sqrt{1+2 \sqrt{2}};
      \end{cases}
    \end{equation}
    
  By calculating the eigenvalues of $B_{2}$, we get 
  \begin{equation}
      \begin{cases}
          \lambda_ {5}^{\prime} =-\frac{1}{2} \sqrt{\frac{1}{2}(-1-2 \sqrt{2}-\mathrm{i} \sqrt{439+164 \sqrt{2}})},
          \\
          \lambda_ {6}^{\prime}=\frac{1}{2} \sqrt{\frac{1}{2}(-1-2 \sqrt{2}-\mathrm{i} \sqrt{439+164 \sqrt{2}})},
          \\
          \lambda_ {7}^{\prime}=-\frac{1}{2} \sqrt{\frac{1}{2}(-1-2 \sqrt{2}+\mathrm{i} \sqrt{439+164 \sqrt{2}})},
          \\
          \lambda_ {8}^{\prime}=\frac{1}{2} \sqrt{\frac{1}{2}(-1-2 \sqrt{2}+\mathrm{i} \sqrt{439+164 \sqrt{2}})};
     \end{cases}
    \end{equation}
  
  By calculating the eigenvalues of $B_{3}$, we get 
    \begin{equation}
      \begin{cases}
          \lambda_ {9}^{\prime} =\lambda_ {10}^{\prime}=-\mathrm{i} \sqrt{\frac{1}{2}(3+\sqrt{2})},
          \\
          \lambda_ {11}^{\prime}=\lambda_ {12}^{\prime}=\mathrm{i} \sqrt{\frac{1}{2}(3+\sqrt{2})};
     \end{cases}
    \end{equation}
  
  By calculating the eigenvalues of $B_{4}$, we get 
  \begin{equation}
    \begin{cases}
  \lambda_ {13}^{\prime} =-\sqrt{\frac{1}{2}-\mathrm{i} \sqrt{7(2+\sqrt{2})}},
  \\
  \lambda_ {14}^{\prime}=\sqrt{\frac{1}{2}-\mathrm{i} \sqrt{7(2+\sqrt{2})}},
  \\
  \lambda_ {15}^{\prime}=-\sqrt{\frac{1}{2}+\mathrm{i} \sqrt{7(2+\sqrt{2})}},
  \\
  \lambda_ {16}^{\prime}=\sqrt{\frac{1}{2}+\mathrm{i} \sqrt{7(2+\sqrt{2})}}.
   \end{cases}
  \end{equation}
  
  By analysing the above 16 eigenvalues, we find that 
  $\lambda_ {1}^{\prime},\lambda_ {2}^{\prime}$ are zero roots;

  The eigenvalues $ \lambda_ {5}^{\prime},\lambda_ {6}^{\prime},\lambda_ {7}^{\prime},\lambda_ {8}^{\prime},\lambda_ {13}^{\prime},\lambda_ {14}^{\prime},\lambda_ {15}^{\prime},\lambda_ {16}^{\prime}$
  are imaginary roots with nonzero real parts;
  
  The eigenvalues $\lambda_ {3}^{\prime},\lambda_ {4}^{\prime},\lambda_ {9}^{\prime},\lambda_ {10}^{\prime},\lambda_ {11}^{\prime},\lambda_ {12}^{\prime}$ 
  are pure imaginary characteristic roots;
  hence the corresponding solutions are unstable.
  
  In conclusion, the relative equilibria of the square configuration with equal masses under the Manev potential is unstable.

  
  \subsection{The stability of square problem under the Schwarzschild potential with equal mass}\label{sec: the stab of square Schw}

  In this section, we consider the stability with equal mass of square problem in $\mathbb{R}^2$, $U$ satisfies the equations as follows
      \begin{equation}
        U=\sum\limits_{1 \leq i<j \leq 4} \frac{1}{\left [  \left(x_{i}-x_{j}\right)^{2}+\left(y_{i}-y_{j}\right)^{2}\right ]^{\frac{1 }{2} } }+\frac{1}{\left [  \left(x_{i}-x_{j}\right)^{2}+\left(y_{i}-y_{j}\right)^{2}\right ]^{\frac{3 }{2} } }.  
      \end{equation}
  
  Let the center of mass be at the origin, and
  $ z \in \mathbb{R}^{8}$ be $z=\left(x_{1}, y_{1}, x_{2}, y_{2}, x_{3}, y_{3}, x_{4}, y_{4}\right)^{T},$ 
  we choose the vertices of square $$z_{0}=(1,0,0, 1,-1,0,0,-1)^{T}$$ to solve the equation $\nabla (U+\omega ^2I)=0$,
  we obtain $\omega^2=\frac{7}{16}+\frac{5}{2 \sqrt{2}}$.
  
  At the point $z_{0}$, the Hessian matrix of $H_{6}$ at $z_{0}$ is  
  $$
  H_{6}=
  \begin{pmatrix}
      \begin{smallmatrix}
        \frac{1}{8}(5+11 \sqrt{2}) & 0 & -\frac{11}{8\sqrt{2}} &\frac{21}{8\sqrt{2}} & -\frac{5}{8} & 0 & -\frac{11}{8\sqrt{2}} &-\frac{21}{8\sqrt{2}} \\
        0 & \frac{1}{32}(-7+44 \sqrt{2}) &\frac{21}{8\sqrt{2}} & -\frac{11}{8\sqrt{2}} & 0 & \frac{7}{32} & -\frac{21}{8\sqrt{2}}  & -\frac{11}{8\sqrt{2}}\\
        -\frac{11}{8\sqrt{2}} &\frac{21}{8\sqrt{2}} & \frac{1}{32}(-7+44 \sqrt{2}) & 0 & -\frac{11}{8\sqrt{2}} & -\frac{21}{8\sqrt{2}} & \frac{7}{32} & 0\\
      \frac{21}{8\sqrt{2}} & -\frac{11}{8\sqrt{2}} & 0 & \frac{1}{8}(5+11 \sqrt{2}) &-\frac{21}{8\sqrt{2}} & -\frac{11}{8\sqrt{2}} & 0 & -\frac{5}{8}\\
        -\frac{5}{8} & 0 & -\frac{11}{8\sqrt{2}} &-\frac{21}{8\sqrt{2}} & \frac{1}{8}(5+11 \sqrt{2}) & 0 & -\frac{11}{8\sqrt{2}} &\frac{21}{8\sqrt{2}} \\
        0 & \frac{7}{32} &-\frac{21}{8\sqrt{2}} & -\frac{11}{8\sqrt{2}} & 0 & \frac{1}{32}(-7+44 \sqrt{2}) &\frac{21}{8\sqrt{2}} & -\frac{11}{8\sqrt{2}}\\
        -\frac{11}{8\sqrt{2}} &-\frac{21}{8\sqrt{2}} & \frac{7}{32} & 0 & -\frac{11}{8\sqrt{2}} &\frac{21}{8\sqrt{2}} & \frac{1}{32}(-7+44 \sqrt{2}) & 0 \\
      -\frac{21}{8\sqrt{2}} & -\frac{11}{8\sqrt{2}} & 0 & -\frac{5}{8} &\frac{21}{8\sqrt{2}} & -\frac{11}{8\sqrt{2}} & 0 &\frac{1}{8}(5+11 \sqrt{2})
            \end{smallmatrix}
  \end{pmatrix}.
  $$

As before, we calculate the traces of $H_{6}\mathscr{D}(A)$ for different irreducible group representations to solve the eigenvalues of matrix $H_{6}$ as follows.
  
  \begin{equation}
   \begin{split}
     \begin{cases}
      \lambda_ {1}=\frac{5+16\sqrt{2}}{4} ,
            \\
            \lambda_ {2}=\frac{1}{16}(-7-20 \sqrt{2}),
            \\
            \lambda_ {3}=\frac{1}{4}(5-5 \sqrt{2}),
            \\
            \lambda_ {4}=\frac{1}{16}(-7+64\sqrt{2}),
            \\
            \lambda_ {5}=0,
            \\
            \lambda_ {6}=\frac{11 \sqrt{2}}{4}.
      \end{cases}
    \end{split}
  \end{equation}
  
  Then linearize the second-order equation 
  \begin{equation}\label{equ:linearied second equation5}
    \ddot{x}_{i}=2 \omega J \dot{x_{i}}+\omega^{2} x_{i}+\nabla_{x_{i}} U, x\in \mathbb{R}^2.
  \end{equation} 
  In the new coordinate $E_{\lambda}=\left \{ V_{i}, V_{j}\right \}$, let $y_{i}=(x_{i},\dot{x_{i}})$, then the linearized second order equation are
  $\dot{y_{i}}=B_{i}y_{i},\,i=1,2,3,4$.
  
  The $B_{i}, i=1,2,3,4$ satisfy
  \begin{equation}
  B_{i}=\left(\begin{array}{cc}
    0 & I_{2} \\
    \omega^{2} I_{2}+\frac{1}{m_{i}} \Lambda_{i} & 2 \omega J
    \end{array}\right),i=1,2,3,4.
  \end{equation}
  The $\Lambda_{i}, i=1,2,3,4$ are 
  $$
  \begin{array}{c}
    \Lambda_{1}=\left(\begin{array}{cc}
    \lambda_{1} & 0 \\
    0 & \lambda_{2}
    \end{array}\right), \quad \Lambda_{2}=\left(\begin{array}{cc}
    \lambda_{3} & 0 \\
    0 & \lambda_{4}
    \end{array}\right), \quad 
    \Lambda_{3}=\left(\begin{array}{cc}
    \lambda_{5} & 0 \\
    0 & \lambda_{5}
    \end{array}\right), \quad 
    \Lambda_{4}=\left(\begin{array}{cc}
    \lambda_{6} & 0 \\
    0 & \lambda_{6}
    \end{array}\right)
  \end{array}.
  $$
  By calculating the eigenvalues of matrix $B_{1}$, we get 
  \begin{equation}
      \begin{cases}
          \lambda_ {1}^{\prime} =0,
          \\
          \lambda_ {2}^{\prime}=0,
          \\
          \lambda_ {3}^{\prime}=-\frac{1}{4} \sqrt{1+4 \sqrt{2}},
          \\
          \lambda_ {4}^{\prime}=\frac{1}{4} \sqrt{1+4 \sqrt{2}};
      \end{cases}
    \end{equation}
    
  By calculating the eigenvalues of $B_{2}$, we get 
  \begin{equation}
      \begin{cases}
          \lambda_ {5}^{\prime} =-\frac{1}{4} \sqrt{\frac{1}{2}(-1+4 \sqrt{2}-\mathrm{i}  \sqrt{-33+9080 \sqrt{2}})},
          \\
          \lambda_ {6}^{\prime}=\frac{1}{4} \sqrt{\frac{1}{2}(-1+4 \sqrt{2}-\mathrm{i}  \sqrt{-33+9080 \sqrt{2}})},
          \\
          \lambda_ {7}^{\prime}=-\frac{1}{4} \sqrt{\frac{1}{2}(-1+4 \sqrt{2}+\mathrm{i}  \sqrt{-33+9080 \sqrt{2}})},
          \\
          \lambda_ {8}^{\prime}=\frac{1}{4} \sqrt{\frac{1}{2}(-1+4 \sqrt{2}+\mathrm{i}  \sqrt{-33+9080 \sqrt{2}})};
     \end{cases}
    \end{equation}
  
  By calculating the eigenvalues of $B_{3}$, we get 
    \begin{equation}
      \begin{cases}
          \lambda_ {9}^{\prime} =\lambda_ {10}^{\prime}=-\frac{1}{4} \mathrm{i} \sqrt{23+12 \sqrt{2}},
          \\
          \lambda_ {11}^{\prime}=\lambda_ {12}^{\prime}=\frac{1}{4} \mathrm{i} \sqrt{23+12 \sqrt{2}};
     \end{cases}
    \end{equation}
  
  By calculating the eigenvalues of $B_{4}$, we get 
  \begin{equation}
    \begin{cases}
      \lambda_ {13}^{\prime} =-\sqrt{-\frac{7}{16}+\frac{3}{\sqrt{2}}-\frac{1}{4} \mathrm{i} \sqrt{440+77 \sqrt{2}}},
      \\
      \lambda_ {14}^{\prime}=\sqrt{-\frac{7}{16}+\frac{3}{\sqrt{2}}-\frac{1}{4} \mathrm{i} \sqrt{440+77 \sqrt{2}}},
      \\
      \lambda_ {15}^{\prime}=-\sqrt{-\frac{7}{16}+\frac{3}{\sqrt{2}}+\frac{1}{4} \mathrm{i} \sqrt{440+77 \sqrt{2}}},
      \\
      \lambda_ {16}^{\prime}=\sqrt{-\frac{7}{16}+\frac{3}{\sqrt{2}}+\frac{1}{4} \mathrm{i} \sqrt{440+77 \sqrt{2}}}.
   \end{cases}
  \end{equation}
  
  By analysing the above 16 eigenvalues, we find that 
  $\lambda_ {1}^{\prime},\lambda_ {2}^{\prime}$ are zero roots;
  $\lambda_ {3}^{\prime},\lambda_ {4}^{\prime}$ are real roots;
  $ \lambda_ {5}^{\prime},\lambda_ {6}^{\prime},\lambda_ {7}^{\prime},\lambda_ {8}^{\prime},\lambda_ {13}^{\prime},\lambda_ {14}^{\prime},\lambda_ {15}^{\prime},\lambda_ {16}^{\prime}$
  are imaginary roots with nonzero real parts;
  $\lambda_ {9}^{\prime},\lambda_ {10}^{\prime},\lambda_ {11}^{\prime},\lambda_ {12}^{\prime}$ 
  are pure imaginary characteristic roots, 
  hence the corresponding solutions are unstable.
  
  In conclusion, under the Schwarzschild potential, the relative equilibria of the square configuration with equal mass is unstable.


\bibliography{ref.bib}

\begin{thebibliography}{10}

\bibitem{Albouy1996}
{\sc Albouy, A.}
\newblock The symmetric central configurations of four equal masses.
\newblock In {\em Hamiltonian dynamics and celestial mechanics ({S}eattle,
  {WA}, 1995)}, vol.~198 of {\em Contemp. Math.} Amer. Math. Soc., Providence,
  RI, 1996, pp.~131--135.

\bibitem{Albouy2012}
{\sc Albouy, A., and Kaloshin, V.}
\newblock Finiteness of central configurations of five bodies in the plane.
\newblock {\em Ann. of Math. (2) 176}, 1 (2012), 535--588.

\bibitem{Arredondo2014}
{\sc Arredondo, J.~A., P{\'e}rez-Chavela, E., and Stoica, C.}
\newblock Dynamics in the schwarzschild isosceles three body problem.
\newblock {\em Journal of Nonlinear Science 24}, 6 (2014), 997--1032.

\bibitem{Atiyah1975}
{\sc Atiyah, M.~F., Patodi, V.~K., and Singer, I.~M.}
\newblock Spectral asymmetry and {R}iemannian geometry. {I}.
\newblock {\em Math. Proc. Cambridge Philos. Soc. 77\/} (1975), 43--69.

\bibitem{Vivina2014}
{\sc Barutello, V.~L., Jadanza, R.~D., and Portaluri, A.}
\newblock Linear instability of relative equilibria for n-body problems in the
  plane.
\newblock {\em Journal of Differential Equations 257}, 6 (2014), 1773--1813.

\bibitem{Ferrario2008}
{\sc Ferrario, D.~L., and Portaluri, A.}
\newblock On the dihedral {$n$}-body problem.
\newblock {\em Nonlinearity 21}, 6 (2008), 1307--1321.

\bibitem{Hampton2006}
{\sc Hampton, M., and Moeckel, R.}
\newblock Finiteness of relative equilibria of the four-body problem.
\newblock {\em Invent. Math. 163}, 2 (2006), 289--312.

\bibitem{Manev2009}
{\sc Haranas, I.~I., and Mioc, V.}
\newblock Manev potential and satellite orbits.
\newblock {\em Romanian Astronomical Journal 19\/} (2009).

\bibitem{Hu2009}
{\sc Hu, X., and Sun, S.}
\newblock Stability of relative equilibria and {M}orse index of central
  configurations.
\newblock {\em C. R. Math. Acad. Sci. Paris 347}, 21-22 (2009), 1309--1312.

\bibitem{Moeckel1992}
{\sc Moeckel, R.}
\newblock A nonintegrable model in general relativity.
\newblock {\em Comm. Math. Phys. 150}, 2 (1992), 415--430.

\bibitem{Moeckel1995}
{\sc Moeckel, R.}
\newblock Linear stability analysis of some symmetrical classes of relative
  equilibria.
\newblock In {\em Hamiltonian dynamical systems ({C}incinnati, {OH}, 1992)},
  vol.~63 of {\em IMA Vol. Math. Appl.} Springer, New York, 1995, pp.~291--317.

\bibitem{Saari2005}
{\sc Saari, D.~G.}
\newblock {\em Collisions, rings, and other {N}ewtonian {$N$}-body problems},
  vol.~104 of {\em CBMS Regional Conference Series in Mathematics}.
\newblock Published for the Conference Board of the Mathematical Sciences,
  Washington, DC; by the American Mathematical Society, Providence, RI, 2005.

\bibitem{Santoprete2006}
{\sc Santoprete, M.}
\newblock Linear stability of the {L}agrangian triangle solutions for
  quasihomogeneous potentials.
\newblock {\em Celestial Mech. Dynam. Astronom. 94}, 1 (2006), 17--35.

\bibitem{Siegel1971}
{\sc Siegel, C.~L., and Moser, J.~K.}
\newblock {\em Lectures on celestial mechanics}.
\newblock Die Grundlehren der mathematischen Wissenschaften, Band 187.
  Springer-Verlag, New York-Heidelberg, 1971.
\newblock Translation by Charles I. Kalme.

\bibitem{Steinberg2012}
{\sc Steinberg, B.}
\newblock {\em Representation theory of finite groups}.
\newblock Universitext. Springer, New York, 2012.
\newblock An introductory approach.

\bibitem{Xia2008}
{\sc Xia, Z.}
\newblock Symmetries in n-body problem.
\newblock {\em EXPLORING THE SOLAR SYSTEM AND THE UNIVERSE 1043\/} (APR 2008),
  126--132.

\bibitem{Xia2021}
{\sc Xia, Z., and Zhou, T.}
\newblock Applying the symmetry groups to study the n body problem.
\newblock {\em Journal of Differential Equations\/} (2021).

\end{thebibliography}

\bibliographystyle{acm}

\end{document}